\documentclass[a4paper]{article}
\usepackage[english]{babel}
\usepackage{pifont,setspace,hyperref,amsthm,amsmath,amssymb,color,multirow,graphicx,subfigure }
\usepackage{enumitem}
\usepackage{inputenc}
\usepackage{amsmath,latexsym,amsthm}
\usepackage{graphicx}
\usepackage{amssymb}
\newcommand{\R}{\mathbb{R}}

\newcommand{\N}{\mathbb{N}}

\newcommand{\capd}{_cD_{0t}^{\alpha}}

\begin{document}
			
				\begin{center}
					\textbf{\large A Boundary-Value Problem for the 3-D Fractional Wave Equation with Singularity}\\
					
					\medskip
					
					{Joseph David, Alexander Nolte and Julie Sherman}
				\end{center}
				
			\makeatletter
			\renewcommand{\@evenhead}{\vbox{\thepage \hfil {\it David J., Nolte A., Sherman J.}   \hrule }}
			\renewcommand{\@oddhead}{\vbox{\hfill
					{\it A Boundary-Value Problem for the 3-D Fractional Wave Equation with Singularity }\hfill \thepage \hrule}} \makeatother
			%
			%
			\label{firstpage}
			
			%
			%
	\begin{abstract}
	In this paper, a boundary-value problem for the 3-D wave equation with Caputo and Bessel operators is investigated. Sufficient conditions on the initial data are established for the existence of a unique solution to the considered problem.
	\end{abstract}		
			\textbf{1. Introduction.}
			
			The field of fractional calculus is concerned with the investigation of integrals and derivatives of arbitrary orders. The basic concepts had been speculated upon by mathematicians such as Leibniz and Euler, but it was not until about 30 years ago that fractional calculus gained specialized interest. This interest grew as applications of fractional calculus to physics, biology, engineering, and other fields were discovered and developed (see \cite{KilSrivTru} - \cite{JesusMach}). Today, there are numerous books and conferences devoted solely to fractional calculus and its applications (see \cite{MillerRoss} - \cite{Ross}). \medskip
			
			The study of fractional diffusion-wave equations (FDWEs) is a subfield of fractional calculus that arises when one substitutes the first and second-order time derivatives in the diffusion and wave equations, respectively, with a fractional derivative of order $\alpha$, $0<\alpha < 2$. In (see \cite{Mainardi}), Mainardi defined the fractional diffusion and fractional wave equations in two dimensions by $$ \frac{\partial^\alpha u}{\partial t^\alpha} = c \frac{\partial ^2 \alpha}{\partial x^2},$$ where fractional diffusion is when $0<\alpha<1$ and the fractional wave equation is when $1<\alpha<2$, and $c>0$. Note that the fractional derivative $\frac{\partial^\alpha u}{\partial t^\alpha}$ has several definitions, such as those by Caputo, Riemann-Liouville, Weyl, and Riesz (see \cite{MillerRoss}, \cite{Podlubny}, \cite{OldSpan}). FDWEs have been shown to model many processes more accurately than classical equations, and, thus, it is desirable to find analytical and numerical solutions\cite{Hilfer}. 
			
			The current literature contains many papers concerned with solving FDWEs of various forms. Some of these differences include the boundary conditions, whether the problem is being considered in arbitrary dimension, and the definition of fractional derivative used (see \cite{BoyLuch} - \cite{Pov}). It has also been considered when the wave equation slightly modified, such as having a finite linear combination of time-fractional derivatives (see \cite{BazBaz}), or where the constant term is replaced by a function of a given form (see \cite{MartReySosa}).  
			
			We note works \cite{Al-Musalhi}, \cite{Ruzhansky}, where direct and inverse problems for PDE with Caputo and Bessel operators were investigated.
			
			In this paper, we will be consider the following partial differential equation 
			\begin{equation}\label{eq1}
			\capd u(x,y,t) = u_{xx} + \frac{1}{x} u_x + u_{yy} + f(x,y,t)
			\end{equation}
			of fractional order $1<\alpha<2$ defined on $\Omega = \{(x,y,t): 0 \leq x, y \leq 1, \,\, 0 \leq t\leq T \}$ with boundary conditions 
			\begin{align}\label{original-conditions-start}
			\lim_{x\rightarrow 0} xu_x &= 0 &&u(1,y,t) = 0\\
			u(x,0,t) &= 0 &&u(x,1,t) = 0
			\end{align} and initial conditions
			\begin{align}
			u(x,y,0) &= \varphi(x,y) &&u_t(x,y,0) = \psi(x,y).\label{original-conditions-end}
			\end{align}
			(For definitions, see section 2). We will find a formal solution, prove uniqueness, establish sufficient conditions for existence, and prove continuous dependence on initial conditions. A variant of the system we consider with zero initial conditions has appeared in \cite{OzKar}, though this paper only writes down a formal solution. The generality of our system, proofs of uniqueness, existence, and continuous dependence are original. Furthermore, this paper is purposefully written with an attention to detail, aiming to be accessible to mathematicians interested in fractional derivatives but without an in-depth background in fractional calculus.

			The paper is structured as followed. In Section 2, basic definitions and properties of the Mittag-Leffler function and Caputo fractional derivative are given, as well as a proof of the solution to a Cauchy problem that arises in finding a formal solution to our FDWE. In Section 3, a formal solution to the problem obtained using the method of separation of variables. In Section 4, the uniqueness of the solution is proved. In Section 5, sufficient conditions for the existence of the solution are determined. Lastly, in Section 6, continuous dependence of our solution on the initial conditions is shown.

			\textbf{2. Preliminaries}

			\textbf{2.1. Basic Definitions}
			
			The Caputo fractional derivative operator of order $\alpha$, where $n\in \N$ and $n-1<\alpha<n$, is defined by 
			\begin{equation*}
			_cD^{\alpha}_{at} g(x) = \frac{1}{\Gamma(n-\alpha)}\int_a^x \frac{g^{(n)}(z)dz}{(x-z)^{\alpha+1-n}}.
			\end{equation*} 
			In this paper we will be working with the specific case of $1<\alpha<2$ and $a=0$, giving
			\begin{equation}\label{caputo-def}
			\capd g(x) = \frac{1}{\Gamma(2-\alpha)}\int_0^x \frac{g''(z)dz}{(x-z)^{\alpha-1}}.
			\end{equation} 
			It has the following relationship with the Reimann-Liouville fractional derivative (see \cite{Podlubny} eq. 2.166)
			\begin{equation}\label{eqncaptorld}
			_{c}D_{ax}^{\alpha} y(x)= _{RL}D_{ax}^{\alpha}y(x) - \sum_{k = 0}^{n-1}\frac{y^{(k)}(a)}{\Gamma(k-\alpha + 1)}(x-a)^{k-\alpha}
			\end{equation}
			where
			\begin{equation*}
						_{RL}D^{\alpha}_{at} g(x) = \frac{1}{\Gamma(n-\alpha)}\frac{d^n}{dx^n}\int_a^x \frac{g(z)dz}{(x-z)^{\alpha+1-n}}
						\end{equation*}
			is the Riemann-Liouville fractional derivative operator of order $\alpha$, where $n\in \N$ and $n-1<\alpha<n$.
			
			An important function in the solution of many differential equations involving fractional derivatives is the Mittag-Leffler function, a generalization of the exponential function. It is defined for all $\alpha > 0$ and $\beta \in \mathbb{R}$ by 
			\begin{equation*}
			E_{\alpha,\beta}(z)=\sum_{k=0}^{\infty} \frac{z^k}{\Gamma(\alpha k + \beta)},
			\end{equation*}
			and has the useful properties that
			\begin{equation}\label{eqnrld}
			_{RL}D_{0t}^{\gamma}t^{\alpha k + \beta - 1}E_{\alpha, \beta}^{(k)}(\lambda t^{\alpha})  = t^{\alpha k + \beta -\gamma- 1}E_{\alpha, \beta-\gamma}^{(k)}(\lambda t^{\alpha}), \hspace{1cm} (\gamma, \beta \in \mathbb{R}, \ \alpha > 0)
			\end{equation}
			(\cite{Podlubny}, eq. 1.82), and for $\alpha > 0$,
			\begin{equation} \label{MittagDeriv}
			\frac{d}{dx}E_{\alpha, 1}(\lambda x^{\alpha}) = \lambda x^{\alpha-1}E_{\alpha, \alpha}(\lambda x^{\alpha}).
			\end{equation}
			
			Here $E_{\alpha, \beta}^{(k)}(\lambda t^{\alpha})=\frac{d^k}{dt^k}E_{\alpha, \beta}(\lambda t^{\alpha})$.
			
			More generally, whenever $k \neq 1$
			\begin{equation} \label{MittagDerivGeneral}
			\frac{d}{dz} z^{k-1}E_{\alpha,k}(\lambda z^\alpha) = z^{k-2}E_{\alpha,k-1}(\lambda z^\alpha). 
			\end{equation} Here, (\ref{MittagDeriv}) - (\ref{godsend-identity}) follow nicely from direct calculations.
			Further, for all $\alpha > 0$ and $\beta \in \mathbb{R}$
			\begin{equation}\label{godsend-identity}
			E_{\alpha, \beta}(z) - z E_{\alpha, \alpha + \beta}(z) = \frac{1}{\Gamma(\beta)},
			\end{equation} and for all $\nu, \beta > 0$,
			\begin{equation}\label{gammaint}
			\frac{1}{\Gamma(\nu)}\int_0^z(z-t)^{\nu-1}E_{\alpha, \beta}(\lambda t^\alpha)t^{\beta-1}dt = z^{\beta+\nu-1} E_{\alpha, \beta+\nu}(\lambda z^{\alpha})
			\end{equation} 
			(\cite{Podlubny}, eq. 1.100) where by standard convention $\frac{1}{\Gamma(0)}$ is understood to be $0$. Finally, the Mittag-Leffler function has the following bound, given in \cite{Podlubny} eq. 1.148
			\begin{equation}\label{strong-E-bound} 
			E_{\alpha, \beta}(z) \leq \frac{C}{1 + |z|},
			\end{equation} 
			which holds for $\alpha<2$, all real $\beta$, and for all $z \in \mathbb{C}$ such that $|\arg(z)|\leq \mu$ for some $\mu$ with $\pi\alpha/2<\mu<\min\{\pi,\pi\alpha\}$. Note, in particular, that this holds for all real $z$ for any $\alpha < 2$.
			
			\par We will also come across the Bessel functions in our paper. For all $p \in \mathbb{Z}$, the Bessel function of the first kind of order $p$ is defined by the series
			\begin{equation}\label{bessel-def-first-kind}
			J_p(x) = \sum_{m = 0}^\infty \frac{(-1)^m}{\Gamma(m + 1) \Gamma(p + m + 1)} \left(\frac {x}{2}\right)^{p + 2m}.
			\end{equation} 
			Also, the Bessel function of the second kind of the $0$-th order is given by the series
			\begin{equation}\label{bessel-def-second-kind}
			Y_0(x) = \frac{2}{\pi}J_0(x)\left(\ln \frac{x}{2} + C\right) - \frac{2}{\pi} \sum_{m =1}^{\infty} \left[\frac{(-1)^m}{(m!)^2}\left( \frac{x}{2}\right)^{2m}\sum_{k = 0}^m \frac{1}{m}\right],
			\end{equation} 
			where $C > 0$ is Euler's constant \cite{Tolstov}. We will make use of a number of special properties of the Bessel functions. First, direct calculations show that $J_0(0) = 1$ and $J_n(0)$ is finite for all $n \in \mathbb{Z}$. Additionally, $J_0$ is known to be bounded (see, for instance, \cite{Tolstov}, Ch. 8, Section 9). By \cite{Tolstov}, eq. 7.9-7.10, for all $p \in \mathbb{R}$, \begin{align}\label{bessel-derivative-fact} J_{p-1}(x) - J_{p+1}(x) &= 2\frac{d}{dx}J_p(x), \\ J_{p-1}(x) + J_{p+1}(x) &= \frac{2p}{x}J_p(x). \label{bessel-division-killer} \end{align} Additionally, for any $p \in \mathbb{R}$,
			\begin{align}
			\frac{d}{dx} \left(x^p J_p(x)\right) &= x^p J_{p-1}(x), \label{Bessel-derivative-neg-shift} \\
			\frac{d}{dx} \left(x^{-p}J_p(x) \right) &= -x^{-p} J_{p+1}(x), \label{Bessel-derivative-pos-shift}
			\end{align} which are \cite{Tolstov}, eq. 7.1-7.2. Facts about the zeros of the Bessel function and bounds that can be placed on its integrals will become relevant in section $5$, and will be given as needed.
			\medskip

			
			\textbf{2.2. Cauchy Problem}
			
			The \textit{Cauchy Problem} given by the fractional differential equation
			\begin{equation}\label{cauchy-problem}
			\capd y(x) - \lambda y(x) = f(x), \hspace{1cm} x \in[0,1],\,\, \lambda \in \mathbb{R}
			\end{equation} with initial conditions 
			\begin{align}\label{cauchy-conditions}
			y(0) &= b_0 &&y'(0) = b_1 
			\end{align} appears repeatedly when solving time-related systems in the problem considered in this paper. Our aim is to show that for $1 < \alpha < 2$,
			\begin{equation}\label{cauchysoln}
			y(x) = b_0E_{\alpha, 1}(\lambda x^{\alpha})+b_1xE_{\alpha, 2}(\lambda x^{\alpha})+ \int_0^x (x-t)^{\alpha-1}E_{\alpha, \alpha}[\lambda(x-t)^{\alpha}]f(t)dt
			\end{equation}
			is a solution to (\ref{cauchy-problem}) which satisfies (\ref{cauchy-conditions}). The following solution is also available, for example, in \cite{KilSrivTru}. 
			Substituting (\ref{cauchysoln}) to the fractional derivative yields 
			\begin{equation}\label{cauchysubsoln} 
			_cD_{0t}^{\alpha} y(x)  =\, _cD_{0t}^{\alpha}\left( b_0E_{\alpha, 1}(\lambda x^{\alpha})\right)+_cD_{0t}^{\alpha}\left(b_1xE_{\alpha, 2}(\lambda x^{\alpha})\right)+ _cD_{0t}^{\alpha}\left(\int_0^x (x-t)^{\alpha-1}E_{\alpha, \alpha}[\lambda(x-t)^{\alpha}]f(t)dt\right).
			\end{equation}
			Consider the first term of (\ref{cauchysubsoln}). We will show that $\capd (b_0 E_{\alpha, 1}(\lambda x^{\alpha})) - \lambda b_0 E_{\alpha, 1}(\lambda x^{\alpha})=0$. Using relation (\ref{eqncaptorld}) with $n = 2$ and $a = 0$ demonstrates that
			\begin{align*}
			\capd (b_0 E_{\alpha, 1}(\lambda x^{\alpha})) 
			& =b_0 \left(_{RL}D_{0t}^{\alpha}(E_{\alpha, 1}(\lambda x^{\alpha})) - \frac{E_{\alpha, 1}(\lambda x^{\alpha})|_{x = 0}}{\Gamma(1-\alpha)}x^{-\alpha} -\frac{\frac{d}{dx}E_{\alpha, 1}(\lambda x^{\alpha})|_{x = 0}}{\Gamma(2-\alpha)}x^{1-\alpha}\right).
			\end{align*}
			From property (\ref{eqnrld}) with $k = 0$, $\gamma = \alpha$, and $\beta = 1$ we have
			$$_{RL}D_{0t}^{\alpha}(E_{\alpha, 1}(\lambda x^{\alpha}))  =  x^{ -\alpha}E_{\alpha, 1-\alpha}(\lambda x^{\alpha}).
			$$
			Thus, 
			\begin{align*}
			\capd (b_0 E_{\alpha, 1}(\lambda x^{\alpha})) 
			& =b_0 \left(x^{ -\alpha}E_{\alpha, 1-\alpha}(\lambda x^{\alpha}) - \frac{x^{-\alpha}}{\Gamma(1-\alpha)} -\frac{\frac{d}{dx}E_{\alpha, 1}(\lambda x^{\alpha})|_{x = 0}}{\Gamma(2-\alpha)}x^{1-\alpha}\right).
			\end{align*}
			Evaluating at $x = 0$ gives $\frac{d}{dx} E_{\alpha, 1}(\lambda x^{\alpha})|_{x=0} = 0 $ because $\alpha >1$. Therefore,
			\begin{align*}
			\capd (b_0 E_{\alpha, 1}(\lambda x^{\alpha})) 
			& =b_0 x^{-\alpha}\left(E_{\alpha, 1-\alpha}(\lambda x^{\alpha}) - \frac{1}{\Gamma(1-\alpha)}\right)\\
			& =b_0 x^{-\alpha}\left(\sum_{k = 0}^{\infty}\frac{(\lambda x^{\alpha})^k}{\Gamma(k\alpha+1-\alpha)}- \frac{1}{\Gamma(1-\alpha)}\right)\\
			& =b_0 \left(\frac{x^{-\alpha}}{\Gamma(1-\alpha)}+\frac{\lambda }{\Gamma(1)}+\frac{\lambda^2 x^{\alpha}}{\Gamma(\alpha+1)}+ \dots- \frac{x^{-\alpha}}{\Gamma(1-\alpha)}\right)\\
			& =b_0 \lambda \left(\frac{1}{\Gamma(1)}+\frac{\lambda x^{\alpha}}{\Gamma(\alpha+1)}+ \dots\right)\\
			& =b_0 \lambda \left(\sum_{k = 0}^{\infty}\frac{(\lambda x^{\alpha})^k}{\Gamma(k\alpha+1)}\right)\\
			& =b_0 \lambda E_{\alpha, 1}(\lambda x^{\alpha})
			\end{align*}
			Indeed, $\capd (b_0 E_{\alpha, 1}(\lambda x^{\alpha})) - \lambda b_0 E_{\alpha, 1}(\lambda x^{\alpha}) = 0$, so the problem is reduced to showing that 
			\begin{equation*}
			b_1xE_{\alpha, 2}(\lambda x^{\alpha})+ \int_0^x (x-t)^{\alpha-1}E_{\alpha, \alpha}[\lambda(x-t)^{\alpha}]f(t)dt
			\end{equation*} is a solution to (\ref{cauchy-problem}).

			
			Consider $_cD_{0t}^{\alpha}(b_1xE_{\alpha, 2}(\lambda x^{\alpha}) )$.  Again, we will show that $\capd (b_1xE_{\alpha, 2}(\lambda x^{\alpha})) - \lambda b_1xE_{\alpha, 2}(\lambda x^{\alpha})=0$. Use relation (\ref{eqncaptorld}) with $n = 2$ and $a = 0$ to see \begin{align*}
			_cD_{0t}^{\alpha}(b_1xE_{\alpha, 2}(\lambda x^{\alpha}) )& =b_1 \left(_{RL}D_{0t}^{\alpha}(xE_{\alpha, 2}(\lambda x^{\alpha})) - \frac{xE_{\alpha, 2}(\lambda x^{\alpha})|_{x = 0}}{\Gamma(1-\alpha)}x^{-\alpha} -\frac{\frac{d}{dx}(xE_{\alpha, 2}(\lambda x^{\alpha}))|_{x = 0}}{\Gamma(2-\alpha)}x^{1-\alpha}\right)\\
			& =b_1 \left(_{RL}D_{0t}^{\alpha}(xE_{\alpha, 2}(\lambda x^{\alpha}))  -\frac{E_{\alpha, 2}(\lambda x^{\alpha})|_{x = 0}}{\Gamma(2-\alpha)}x^{1-\alpha}\right)\\
			& =b_1 \left(_{RL}D_{0t}^{\alpha}(xE_{\alpha, 2}(\lambda x^{\alpha}))  -\frac{x^{1-\alpha}}{\Gamma(2-\alpha)}\right).
			\end{align*}
			Further, by property (\ref{eqnrld}) with $k = 0$, $\gamma = \alpha$, and $\beta = 2$,
			$$
			_{RL}D_{0t}^{\alpha}(xE_{\alpha, 2}(\lambda x^{\alpha}))  =  x^{ 1-\alpha}E_{\alpha, 2-\alpha}(\lambda x^{\alpha}).
			$$
			Thus, 	
			\begin{align*}
			_cD_{0t}^{\alpha}(b_1xE_{\alpha, 2}(\lambda x^{\alpha}) )
			& =b_1 x^{ 1-\alpha}\left(E_{\alpha, 2-\alpha}(\lambda x^{\alpha})  -\frac{1}{\Gamma(2-\alpha)}\right)\\
			& =b_1 x^{1-\alpha}\left(\sum_{k = 0}^{\infty}\frac{(\lambda x^{\alpha})^k}{\Gamma(k\alpha+2-\alpha)}- \frac{1}{\Gamma(2-\alpha)}\right)\\
			& =b_1 x^{1-\alpha}\left(\frac{1}{\Gamma(2-\alpha)}+\frac{\lambda x^{\alpha}}{\Gamma(2)}+\frac{\lambda^2 x^{2\alpha}}{\Gamma(\alpha+2)}+ \dots- \frac{1}{\Gamma(2-\alpha)}\right)\\
			& =b_1 \lambda x\left(\frac{1}{\Gamma(2)}+\frac{\lambda x^{\alpha}}{\Gamma(\alpha+2)}+ \dots\right)\\
			& =b_1 \lambda x \left(\sum_{k = 0}^{\infty}\frac{(\lambda x^{\alpha})^k}{\Gamma(k\alpha+2)}\right)\\
			& =b_1 \lambda xE_{\alpha, 2}(\lambda x^{\alpha})
			\end{align*}
			Indeed, $\capd b_1xE_{\alpha, 2}(\lambda x^{\alpha}) - \lambda b_1xE_{\alpha, 2}(\lambda x^{\alpha}) = 0$, so the problem is further reduced to showing that 
			\begin{equation*}
			\int_0^x (x-t)^{\alpha-1}E_{\alpha, \alpha}[\lambda(x-t)^{\alpha}]f(t)dt
			\end{equation*} is a solution to (\ref{cauchy-problem}).

			

			By the definition (\ref{caputo-def}),
			\begin{equation}\label{caputoint}
			\begin{array}{l}
			\displaystyle{_cD_{0t}^{\alpha}\left(\int_0^x (x-t)^{\alpha-1}E_{\alpha, \alpha}[\lambda(x-t)^{\alpha}]f(t)dt\right)  =}\\ \displaystyle{\frac{1}{\Gamma(2-\alpha)}\int_0^x\frac{1}{(x-t)^{\alpha-1}}\frac{d^2}{dt^2}\left(\int_0^t (t-z)^{\alpha-1}E_{\alpha, \alpha}[\lambda(t-z)^{\alpha}]f(z)dz\right)dt.}
			\end{array}
			\end{equation}
			We will need the following well-known equation,
			\begin{equation}\label{dirchint}
			\frac{d}{dx}\int_0^xf(x,t)dt = f(x,x) + \int_0^x\frac{d}{dx}f(x,t)dt. 
			\end{equation}
			We proceed by evaluating the derivative that appears in (\ref{caputoint}). By (\ref{dirchint}),
			\begin{align*}
			&\frac{d}{dt}\left(\int_0^t (t-z)^{\alpha-1}E_{\alpha, \alpha}[\lambda(t-z)^{\alpha}]f(z)dz\right) = (t-t)^{\alpha-1}E_{\alpha, \alpha}[\lambda(t-t)^{\alpha}]f(t) + \\
			&+\int_0^t \frac{d}{dt}\left((t-z)^{\alpha-1}E_{\alpha, \alpha}[\lambda(t-z)^{\alpha}]f(z)\right)dz
			=  \int_0^t \frac{d}{dt}\left((t-z)^{\alpha-1}E_{\alpha, \alpha}[\lambda(t-z)^{\alpha}]f(z)\right)dz.
			\end{align*}
			To simplify the integrand on the right hand side, use (\ref{eqnrld}) with $\gamma = 1$, $k = 0$, and $\beta = \alpha$ to see
			\begin{align*}
			\frac{d}{dt}\left(\int_0^t (t-z)^{\alpha-1}E_{\alpha, \alpha}[\lambda(t-z)^{\alpha}]f(z)dz\right)
			& =  \int_0^t (t-z)^{\alpha - 2}E_{\alpha, \alpha-1}(\lambda (t-z)^{\alpha})f(z)dz.
			\end{align*}
			
			Now we want to consider the second derivative
			\begin{align*}
			\frac{d^2}{dt^2}\left(\int_0^t (t-z)^{\alpha-1}E_{\alpha, \alpha}[\lambda(t-z)^{\alpha}]f(z)dz\right) & = \frac{d}{dt}\int_0^t (t-z)^{\alpha - 2}E_{\alpha, \alpha-1}(\lambda (t-z)^{\alpha})f(z)dz.
			\end{align*} By changing variables so that $z = t-z$ and applying (\ref{dirchint}),
			\begin{align}
			\frac{d^2}{dt^2}\left(\int_0^t (t-z)^{\alpha-1}E_{\alpha, \alpha}[\lambda(t-z)^{\alpha}]f(z)dz\right) & = \frac{d}{dt}\int_0^t z^{\alpha - 2}E_{\alpha, \alpha-1}(\lambda z^{\alpha})f(t-z)dz\\
			& = t^{\alpha - 2}E_{\alpha, \alpha-1}(\lambda t^{\alpha})f(t-t)+\int_0^t \frac{d}{dt}\left(z^{\alpha - 2}E_{\alpha, \alpha-1}(\lambda z^{\alpha})f(t-z)\right)dz\\
			& = t^{\alpha - 2}E_{\alpha, \alpha-1}(\lambda t^{\alpha})f(0)+\int_0^t \frac{d}{dt}\left(z^{\alpha - 2}E_{\alpha, \alpha-1}(\lambda z^{\alpha})f(t-z)\right)dz.\label{caputo-second-der-end}
			\end{align}
			Now, plugging back into equation (\ref{caputoint})
			\begin{align}
			_c&D_{0t}^{\alpha}\left(\int_0^x (x-t)^{\alpha-1}E_{\alpha, \alpha}[\lambda(x-t)^{\alpha}]f(t)dt\right) \\
			& =\frac{1}{\Gamma(2-\alpha)}\int_0^x\frac{1}{(x-t)^{\alpha-1}}\left(t^{\alpha - 2}E_{\alpha, \alpha-1}(\lambda t^{\alpha})f(0)+\int_0^t \frac{d}{dt}\left(z^{\alpha - 2}E_{\alpha, \alpha-1}(\lambda z^{\alpha})f(t-z)\right)dz\right)dt\\
			& = \frac{1}{\Gamma(2-\alpha)}\left[\int_0^x\frac{1}{(x-t)^{\alpha-1}}t^{\alpha - 2}E_{\alpha, \alpha-1}(\lambda t^{\alpha})f(0)dt+
			\int_0^x\frac{1}{(x-t)^{\alpha-1}}\int_0^t \frac{d}{dt}\left(z^{\alpha - 2}E_{\alpha, \alpha-1}(\lambda z^{\alpha})f(t-z)\right)dz dt\right] \label{yucky}
			\end{align}

			Consider the first term of (\ref{yucky}). By (\ref{gammaint}) with $\nu = 2 - \alpha$ and $\beta = \alpha - 1$, we obtain
			\begin{align*}
			\frac{1}{\Gamma(2-\alpha)}\int_0^x\frac{1}{(x-t)^{\alpha-1}}t^{\alpha - 2}E_{\alpha, \alpha-1}(\lambda t^{\alpha})f(0)dt & = f(0)E_{\alpha, 1}(\lambda x^{\alpha}).
			\end{align*}
			Finally, consider the second term of (\ref{yucky}).
			In the inner integral, change variables by $z = t-z$ to see
			\begin{align*}
			&\frac{1}{\Gamma(2-\alpha)}\int_0^x\frac{1}{(x-t)^{\alpha-1}}\int_0^t z^{\alpha - 2} E_{\alpha, \alpha-1}(\lambda z^{\alpha})\frac{d}{dt}f(t-z)dz dt\\
			& = \frac{1}{\Gamma(2-\alpha)}\int_0^x{(x-t)^{1-\alpha}}\int_0^t (t-z)^{\alpha - 2}E_{\alpha, \alpha-1}(\lambda (t-z)^{\alpha})f'(z)dz dt
			\end{align*}
			Now by Fubini's theorem,
			\begin{align*}
			&\frac{1}{\Gamma(2-\alpha)}\int_0^x{(x-t)^{1-\alpha}}\int_0^t (t-z)^{\alpha - 2} E_{\alpha, \alpha-1}(\lambda (t-z)^{\alpha})f'(z)dz dt\\
			& = \frac{1}{\Gamma(2-\alpha)}\int_0^xf'(z)\int_z^x {(x-t)^{1-\alpha}}(t-z)^{\alpha - 2}E_{\alpha, \alpha-1}(\lambda (t-z)^{\alpha}) dtdz
			\end{align*}
			In the inner integral, setting $\xi = t-z$ yields
			\begin{align*}
			&\frac{1}{\Gamma(2-\alpha)}\int_z^x {(x-t)^{1-\alpha}}(t-z)^{\alpha - 2}E_{\alpha, \alpha-1}(\lambda (t-z)^{\alpha}) dt \\
			& = \frac{1}{\Gamma(2-\alpha)}\int_0^{x-z} {(x-(\xi + z))^{1-\alpha}}\xi^{\alpha - 2}E_{\alpha, \alpha-1}(\lambda \xi^{\alpha}) d\xi.
			\end{align*} Using (\ref{gammaint}) again with $\nu = 2 - \alpha$ and $\beta = \alpha - 1$ and $z = x-z$ demonstrates
			\begin{align*}
			\frac{1}{\Gamma(2-\alpha)}\int_z^x {(x-t)^{1-\alpha}}(t-z)^{\alpha - 2}E_{\alpha, \alpha-1}(\lambda (t-z)^{\alpha}) dt & = E_{\alpha, 1} (\lambda  (x-z)^{\alpha})
			\end{align*}
			Thus, from equation (\ref{caputoint}) we now have
			\begin{align*}
			_cD_{0t}^{\alpha}\left(\int_0^x (x-t)^{\alpha-1}E_{\alpha, \alpha}[\lambda(x-t)^{\alpha}]f(t)dt\right) & = f(0)E_{\alpha, 1}(\lambda x^{\alpha}) + \int_0^xf'(z)E_{\alpha, 1} (\lambda  (x-z)^{\alpha})dz.
			\end{align*}
			After using integration by parts on the above integral, we see 
			\begin{align*}
			\int_0^xf'(z)E_{\alpha, 1} (\lambda  (x-z)^{\alpha})dz & = f(z)E_{\alpha, 1} (\lambda  (x-z)^{\alpha})\Big|_0^x-\int_0^xf(z)\frac{d}{dz}[E_{\alpha, 1} (\lambda  (x-z)^{\alpha})]dz\\
			& = f(x)E_{\alpha, 1} (0)-f(0)E_{\alpha, 1} (\lambda  x^{\alpha})-\int_0^xf(z)\frac{d}{dz}[E_{\alpha, 1} (\lambda  (x-z)^{\alpha})]dz\\
			& = f(x)-f(0)E_{\alpha, 1} (\lambda  x^{\alpha})-\int_0^xf(z)\frac{d}{dz}[E_{\alpha, 1} (\lambda  (x-z)^{\alpha})]dz.
			\end{align*}
			
			We have now reduced (\ref{caputoint}) to 
			\begin{align*}
			_cD_{0t}^{\alpha}\left(\int_0^x (x-t)^{\alpha-1}E_{\alpha, \alpha}[\lambda(x-t)^{\alpha}]f(t)dt\right) 
			&= f(x)-\int_0^xf(z)\frac{d}{dz}[E_{\alpha, 1} (\lambda  (x-z)^{\alpha})]dz.
			\end{align*} 
			Thus, it remains only to show that 
			\begin{equation}\label{desired}
			-\int_0^xf(z)\frac{d}{dz}[E_{\alpha, 1} (\lambda  (x-z)^{\alpha})]dz - \lambda\int_0^x (x-z)^{\alpha-1}E_{\alpha, \alpha}[\lambda(x-z)^{\alpha}]f(z)dz = 0
			\end{equation}
			Indeed, by (\ref{MittagDeriv})
			\begin{align*}
			\frac{d}{dz}E_{\alpha, 1} (\lambda  (x-z)^{\alpha}) 
			& = -\lambda (x-z)^{\alpha-1}E_{\alpha, \alpha} (\lambda  (x-z)^{\alpha})
			\end{align*}
			so (\ref{desired}) holds.
			Finally, we have seen that (\ref{cauchysoln}) is a solution to (\ref{cauchy-problem}).
			
			It remains to show that
			\begin{align*}
			y(0) & = b_0, &
			y'(0) & = b_1.
			\end{align*} Certainly
			\begin{align*}
			y(0) &= b_0E_{\alpha, 1}(\lambda\cdot 0^{\alpha})+b_1\cdot 0 \cdot E_{\alpha, 2}(\lambda x^{\alpha})+ \int_0^0 (0-t)^{\alpha-1}E_{\alpha, \alpha}[\lambda(0-t)^{\alpha}]f(t)dt  = b_0.
			\end{align*}
			Further, 
			\begin{align}\label{yprime}
			y'(0) &=\frac{d}{dx} (b_0E_{\alpha, 1}(\lambda x^{\alpha}))|_{x=0}+\frac{d}{dx}(b_1xE_{\alpha, 2}(\lambda x^{\alpha}))|_{x=0}+ \frac{d}{dx}\left(\int_0^x (x-t)^{\alpha-1}E_{\alpha, \alpha}[\lambda(x-t)^{\alpha}]f(t)dt\right)\Big|_{x=0}.
			\end{align}
			Consider the first term of (\ref{yprime}). By (\ref{MittagDeriv}), $$\frac{d}{dx} (b_0E_{\alpha, 1}(\lambda x^{\alpha})) = b_0\lambda x^{\alpha-1}E_{\alpha, \alpha}(\lambda x^{\alpha}).$$ Evaluating at $x = 0$ yields $$\frac{d}{dx} (b_0E_{\alpha, 1}(\lambda x^{\alpha}))|_{x=0} = 0.$$
			Thus, (\ref{yprime}) is reduced to 
			\begin{align}\label{yprime1}
			y'(0) &=\frac{d}{dx}(b_1xE_{\alpha, 2}(\lambda x^{\alpha}))|_{x=0}+ \frac{d}{dx}\left(\int_0^x (x-t)^{\alpha-1}E_{\alpha, \alpha}[\lambda(x-t)^{\alpha}]f(t)dt\right)\Big|_{x=0}
			\end{align}
			Considering the first term of (\ref{yprime1}) shows
			\begin{align*}
			\frac{d}{dx}(b_1xE_{\alpha, 2}(\lambda x^{\alpha}))|_{x=0} & = b_1x\frac{d}{dx}(E_{\alpha, 2}(\lambda x^{\alpha}))|_{x=0} + b_1E_{\alpha, 2}(\lambda x^{\alpha})|_{x=0} = b_1.
			\end{align*}
			Finally, consider the last term of (\ref{yprime1}). By the Fundamental Theorem of Calculus,
			\begin{align*}
			\frac{d}{dx}\left(\int_0^x (x-t)^{\alpha-1}E_{\alpha, \alpha}[\lambda(x-)^{\alpha}]f(t)dt\right)\Big|_{x=0} & =  (x-x)^{\alpha-1}E_{\alpha, \alpha}[\lambda(x-x)^{\alpha}]f(x)\Big|_{x=0} = 0.
			\end{align*}
			Thus we also have that $y'(0) = b_1$. Therefore (\ref{cauchysoln}) is a solution to (\ref{cauchy-problem}) - (\ref{cauchy-conditions}).
			\medskip
			

			\textbf{3. Problem Formulation and Formal Solution}
			
			We will be considering the system given by (\ref{eq1})-(\ref{original-conditions-end}) in this paper.
			
			\par We search for solutions of the form $u(x,y,t) = v(x, y, t) + w(x, y, t)$. Here $v(x,y,t)$ satisfies the \textit{homogeneous problem} given by the differential equation
			\begin{equation}\label{hom-eq}
			\capd v(x,y,t) = v_{xx} + \frac{1}{x} v_x + v_{yy}
			\end{equation}
			with boundary and initial conditions
			\begin{align}\label{hom-conditions-start}
			\lim_{x \rightarrow 0} x v_x &= 0 &&v(1,y,t) = 0 \\
			v(x, 0, t) &= 0  &&v(x, 1, t) = 0\label{hom-y-conds} \\
			v(x, y, 0) &= \varphi(x,y)   &&v_t(x,y,0) = \psi(x,y),\label{hom-conditions-end}
			\end{align} and $w(x,y,t)$ satisfies the \textit{non-homogeneous problem} with homogeneous boundary conditions given by the differential equation 
			\begin{equation}\label{nonhom-eq}
			\capd w(x,y,t) = w_{xx} + \frac{1}{x} w_x + w_{yy} + f(x,y,t)
			\end{equation}
			with boundary and initial conditions
			\begin{align}\label{nonhom-conditions-start}
			\lim_{x \rightarrow 0} x w_x &= 0 &&w(1,y,t) = 0 \\
			w(x, 0, t) &= 0  &&w(x, 1, t) = 0 \\
			w(x, y, 0) &= 0   &&w_t(x,y,0) = 0.\label{nonhom-conditions-end}
			\end{align}
			
			Once we have found such $v(x,y,t)$ and $w(x,y,t)$, then $u(x,y,t) = v(x,y,t) + w(x,y,t)$ will be a solution to the original problem. 
			
			
			\medskip 
			\textbf{3.1. Homogeneous Problem}
			
			Consider the homogeneous problem given by (\ref{hom-eq})-(\ref{hom-conditions-end}). We use the method of separation of variables.
			\par Suppose $$u(x,y,t) = Z(x,y)T(t).$$ Then, $$Z(x,y)\capd T(t) = Z_{xx}(x,y)T(t) + \frac{1}{x} Z_x(x,y) T(t)+  Z_{yy}(x,y)T(t).$$ Using the fact that we are searching only for nonzero solutions, we rearrange to see that \begin{equation}\label{first-separation}\frac{\capd T(t)}{T(t)} = \frac{Z_{xx}(x,y) + \frac{1}{x} Z_x(x,y)+  Z_{yy}(x,y)}{Z(x,y)}.\end{equation} Here, the left-hand side depends only on $t$ and the right-hand side depends only on $x$ and $y$. The only way for these to be equal for all $x, y,$ and $t$ is for each to be constant. We denote such a constant by $-\mu$.
			
			\par We now search for solutions to the spatial problem given by (\ref{first-separation}),
			
			\begin{equation*}
			{Z_{xx}(x,y) + \frac{1}{x} Z_x(x,y)+  Z_{yy}(x,y)+\mu Z(x,y)} = 0.
			\end{equation*} We now apply the method of separation of variables for the spatial problem. Suppose $Z(x,y)$ is of the form $X(x)Y(y).$ Then,
			$$X''(x)Y(y) + \frac{1}{x}X'(x)Y(y) + X(x)Y''(y)+ \mu X(x)Y(y) = 0.$$ Using the fact that we are searching for nonzero solutions, we rearrange to see that \begin{equation}\label{second-separation}\frac{Y''(y)}{Y(y)} = \frac{X''(x)+ \frac{1}{x}X'(x)+\mu X(x)}{-X(x)} = -\lambda\end{equation}
			where $\lambda$ is a constant by the same argument as in (\ref{first-separation}).

			
			Next, consider the function $Y(y)$. Combining (\ref{second-separation}) and (\ref{hom-y-conds}) we see that $$\frac{d^2}{dy^2}Y(y) + \lambda Y(y) = 0$$ and that $Y$ satisfies $$Y(0) = Y(1) = 0.$$ It is well-known and simple to verify that this problem has nontrivial solutions $Y_n(y) = \sin (n \pi y)$ and $\lambda_n = -(n\pi)^2$ where $n$ varies over $\mathbb{N}$.
			
			Now consider the function $X(x)$. From (\ref{hom-conditions-start}) and (\ref{second-separation}) we see that $X$ must satisfy $$\frac{d^2}{dx^2}X(x)+ \frac{1}{x}\frac{d}{dx}X(x)+(\mu -\lambda) X(x) = 0$$ and the boundary condition $$X(1) = \lim_{x\rightarrow 0}x\frac{d}{dx}X(x) = 0.$$ Note that this is Bessel's equation with $p=0$. By \cite{Tolstov}, Ch. 8, Section 5, $X$ must have the form $$X(x) = c_1J_0\left((\mu-\lambda)^{1/2}x\right)+c_2Y_0\left((\mu-\lambda)^{1/2}x\right).$$ Denote $\gamma = (\mu-\lambda)^{1/2}$. We now apply the boundary condition at $0$ with the goal of eliminating the $Y_0$ term. Taking the limit as $x \rightarrow 0$ and applying (\ref{bessel-derivative-fact}) shows that
			\begin{align*}
			\lim_{x\rightarrow 0}x\frac{d}{dx}X(x) &= \lim_{x\rightarrow 0}\left(x \frac{d}{dx}[c_1J_0(\gamma x)+c_2Y_0(\gamma x)]\right)\\
			&= \lim_{x\rightarrow 0}\left[c_1 x\frac{d}{dx}J_0(\gamma x)+c_2 x\frac{d}{dx}Y_0(\gamma x)\right]\\
			&= \lim_{x\rightarrow 0}\left[c_1\gamma x\frac{J_{-1}(\gamma x)-J_{1}(\gamma x)}{2}+c_2 x\frac{d}{dx}Y_0(\gamma x)\right] &&\text{and since $J_{-1}$, $J_1$ are bounded at $0$,} \\
			&= \lim_{x \rightarrow 0} c_2 x \frac{d}{dx}Y_0(\gamma x)
			\end{align*}
			We now use the representation of $Y_0$ given by (\ref{bessel-def-second-kind}). Term-by-term differentiating yields that
			\begin{align*}
			\lim_{x\rightarrow 0}&\left[c_2 x\frac{d}{dx}Y_0(\gamma x)\right]\\
			&= \lim_{x\rightarrow 0}\left[c_2 x \frac{d}{dx}\left(\frac{2}{\pi}J_0(\gamma x) \left(\ln\frac{\gamma x}{2}+C \right) - \frac{2}{\pi} \sum_{m=1}^{\infty} \frac{(-1)^m}{(m!)^2}\left(\frac{\gamma x}{2}\right)^{2m}\left(1+\frac{1}{2} +\frac{1}{3}+ \cdots + \frac{1}{m} \right)\right)\right]\\
			&= \lim_{x\rightarrow 0}\left[\frac{2}{\pi}c_2 x\left(\Big(\frac{d}{dx}J_0(\gamma x) \Big)\left(\ln\frac{\gamma x}{2}+C \right)+\frac{J_0(\gamma x)}{x}  -  \sum_{m=1}^{\infty} \gamma m\frac{(-1)^m}{(m!)^2}\left(\frac{\gamma x}{2}\right)^{2m-1}\left(1+\frac{1}{2} +\frac{1}{3}+ \cdots + \frac{1}{m} \right)\right)\right]\\
			&= \frac{2c_2}{\pi}\lim_{x\rightarrow 0}\left[x\Big(\frac{d}{dx}J_0(\gamma x)\Big)\left(\ln\frac{\gamma x}{2}+C \right) + J_0(\gamma x)  -  \sum_{m=1}^{\infty} \gamma m \frac{(-1)^m}{(m!)^2}\left(\frac{\gamma x}{2}\right)^{2m}\left(1+\frac{1}{2} +\frac{1}{3}+ \cdots + \frac{1}{m} \right)\right].
			\end{align*} Now, since the summation term clearly converges uniformly, we can swap the limit and summation. Since every term of the summand is $0$ at $x = 0$, we obtain
			\begin{align*}
			\lim_{x\rightarrow 0}\left[c_2 x\frac{d}{dx}Y_0(\gamma x)\right]&= \frac{2c_2}{\pi}\lim_{x\rightarrow 0}\left[x\Big(\frac{d}{dx}J_0(\gamma x)\Big)\left(\ln\frac{\gamma x}{2}+C \right) + J_0(\gamma x)  \right]\\
			&= \frac{2c_2}{\pi}\lim_{x\rightarrow 0}\left[x\gamma\frac{J_{-1}(\gamma x)-J_{1}(\gamma x)}{2} \left(\ln\frac{\gamma x}{2}+C \right)+ J_0(\gamma x)  \right]\\
			&= \frac{2c_2}{\pi}\lim_{x\rightarrow 0}J_0(\gamma x).
			\end{align*} The last equality follows from $J_n(0)$ being finite for all integer $n$ and the fact that $\lim_{x\rightarrow 0}(x\ln x)= 0$. Finally, we know that $J_0(0) = 1$, therefore, $$ \displaystyle\lim_{x\rightarrow 0}x\frac{d}{dx}X(x)
			=\frac{2c_2}{\pi}.$$ To satisfy the condition $ \lim_{x\rightarrow 0}xX'(x)
			=0$, we must have $c_2 = 0$. Thus, $X(x)$ must have the form $c_1J_0(\gamma x)$. To satisfy the condition $X(1) = 0$, we must have that $\gamma$ is a zero of $J_0(x)$. Denoting the $m^{th}$ zero of $J_0(x)$ by $\gamma_m$, we obtain the general solution $$X_m(x) = c_mJ_0(\gamma_m x),$$ where $m$ varies over $\mathbb{N}$.\\
			
			
			Finally, we turn to the time problem given by our separation of variables:
			\begin{equation*}
			\capd T(t) + \mu T(t) = 0.
			\end{equation*} Combining the results from the previous two sections we see that the only possible values of $\mu$ are given by $$\mu_{m,n} = (n\pi)^2 + \gamma_m^2$$ where $m$ and $n$ are natural numbers. This is the Cauchy problem discussed in section 2.2 with $\lambda = - \mu_{m,n}$ and $f(x) = 0$. Hence, a general solution for any given $\mu_{m,n}$ is given by $$T_{m,n}(t) = A_{m,n}E_{\alpha, 1}(-\mu_{m,n} t^\alpha) + B_{m,n}t E_{\alpha, 2}(-\mu_{m,n} t^\alpha),$$ where $A_{m,n}$ and $B_{m,n}$ are constants.
			
			
			Bringing together the results from our separations of variables we see that our solution to the homogeneous problem must have the form
			$$ u(x,y,t) = \sum_{m, n \in \N} [C_{m,n} E_{\alpha, 1}(-\mu t^\alpha) + D_{m,n}t E_{\alpha, 2}(-\mu t^\alpha)]J_0(\gamma_m x) \sin(n \pi y)$$ where $C_{m,n}$ and $D_{m,n}$ are constants. We now determine these coefficients by applying our initial conditions. Since the Bessel functions $\{J_0\}$ are orthogonal with weight $x$, we use the $L^2$ inner product with weight $x$ throughout this process and the rest of the paper. Since $x$ is positive, continuous, and only zero at a single point, this is an equivalent inner product to the standard one. At $t = 0$, %
			\begin{align*} 
			\varphi(x, y) &= u(x, y, 0) \\
			&= \sum_{m, n \in \N} [C_{m,n} E_{\alpha, 1}(0) + (0) D_{m,n} E_{\alpha, 2}(0)]J_0(\gamma_m x) \sin(n \pi y) \\
			&= \sum_{m, n \in \N} C_{m,n} J_0(\gamma_m x) \sin(n \pi y).
			\end{align*} The last equality is due to the fact that \begin{equation}\label{E-alpha-one-zero}E_{\alpha, 1}(0) = \sum_{k=1}^\infty \frac{(0 ^ \alpha)^k}{\Gamma(k \alpha + 1)} = \frac{1}{\Gamma(1)} = 1.\end{equation} Now, for any $l, k \in \N$, by taking the $L^2$ inner product with $J_0(\gamma_l x) \sin(k\pi y)$ and using the orthogonality of $\{J_0(\gamma_m x) \sin(n \pi y)\}$ we see that
			\begin{align*}
			\langle \varphi(x,y), J_0(\gamma_l x) \sin(k\pi y) \rangle &= \Big \langle \sum_{m, n \in \N} C_{m,n} J_0(\gamma_m x) \sin(n \pi y), J_0(\gamma_l x) \sin(k\pi y) \Big \rangle \\
			&= \sum_{m, n \in \N} \langle C_{m,n} J_0(\gamma_m x) \sin(n \pi y), J_0(\gamma_l x) \sin(k \pi y) \rangle \\
			&= C_{l, k} \langle J_0(\gamma_l x) \sin(k \pi y), J_0(\gamma_l x) \sin(k \pi y) \rangle \\
			&= C_{l, k} ||J_0(\gamma_l x) \sin(k \pi y)||_2^2,
			\end{align*}
			hence for any $l, k \in \N$, \begin{equation}\label{cmn}C_{l, k} = \frac{\langle \varphi(x,y), J_0(\gamma_l x) \sin(k\pi y) \rangle}{||J_0(\gamma_l x) \sin(k \pi y)||_2^2} = \frac{1}{||J_0(\gamma_l x) \sin(k \pi y)||_2^2} \int_0^1 \int_0^1 x\varphi(x,y) J_0(\gamma_l x)\sin(k \pi y) dx dy.\end{equation}
			
			\par To deduce the values of $D_{m,n}$, we use (\ref{MittagDeriv})
			and note that 
			\begin{equation}
			\label{E-alpha-two-zero}
			E_{\alpha, 2}(0) = \sum_{k=0}^{\infty}\frac{(0^\alpha)^k}{\Gamma(k \alpha + 2)} = \frac{1}{\Gamma(2)} = 1.\end{equation} 
			Then applying our initial condition,
			\begin{align*}
			\psi(x,y) &= \frac{d}{dt} u(x,y,t) \Big|_{t=0} \\
			&= \frac{d}{dt} \bigg[\sum_{m, n \in \N} \big[C_{m,n} E_{\alpha, 1}(-\mu t^\alpha) + D_{m,n}t E_{\alpha, 2}(-\mu t^\alpha)\big]J_0(\gamma_m x) \sin(n \pi y) \bigg]_{t =0} \\
			&= \sum_{m,n}\bigg[ \Big[C_{m,n}\frac{d}{dt} E_{\alpha, 1}(-\mu t^\alpha) + D_{m,n}\frac{d}{dt} t E_{\alpha, 2}(-\mu t^\alpha) \Big]_{t=0} J_0(\gamma_m x) \sin(n \pi y)\bigg] \\
			&= \sum_{m,n} \bigg[\Big[C_{m,n} (-\mu t^{\alpha -1}) E_{\alpha, \alpha}(-\mu t^\alpha) + D_{m,n}\Big(t \frac{d}{dt} E_{\alpha, 2}(-\mu t^\alpha) + E_{\alpha, 2}(-\mu t^\alpha) \Big) \Big]_{t=0} J_0(\gamma_m x) \sin(n \pi y)\bigg] \\
			&= \sum_{m,n} D_{m,n} J_0(\gamma_m x) \sin(n \pi y).
			\end{align*} 
			By identical reasoning to that used when determining the $C_{m,n},$ we see that for any $l, k \in \N$, 
			\begin{equation}\label{dmn}
			D_{l, k} =  \frac{1}{||J_0(\gamma_l x) \sin(k \pi y)||_2^2} \int_0^1 \int_0^1 x\psi(x,y) J_0(\gamma_l x)\sin(k \pi y) dx dy.
			\end{equation} This fully characterizes the solution to the homogeneous problem.\\
			
			
			\textbf{3.2. Non-Homogeneous Problem}
			\par Now we tackle the remaining non-homogeneous problem given by (\ref{nonhom-eq}) and (\ref{nonhom-conditions-start})-(\ref{nonhom-conditions-end}). The first step in this process will be to represent the source function $f(x, y, t)$ in a summation of the form $$\sum_{m,n} C(t) J_0(\gamma_m x) \sin(n \pi y).$$ 
			
			\par To do this, we will need the following well-known lemma. 
			
			\textbf{Lemma 1.} \emph{If $\{f_n\}_{n \in \N}$ and $\{g_m\}_{m \in \N}$ are complete orthonormal families in $L^2([0,1])$ then $\{f_n g_m\}_{m, n \in \N}$ is a complete orthonormal family in $L^2([0,1] \times [0,1])$}
			
			\medskip
			
			{\sf Proof. } Let $\{f_n\}_{n \in \N}$ and $\{g_m\}_{m \in \N}$ be complete orthonormal sets in $L^2([0,1])$. A family of orthonormal elements $\{\phi_\alpha\}_{\alpha \in A}$ of a Hilbert space $H$ is complete if and only if the closure of the set of finite linear combinations of $\{\phi_\alpha\}$ is $H$. Since an element $x \in H$ is in the closure of the set of finite linear combinations of $\{\phi_\alpha\}$ if and only if the Parsevall identity, $$ ||x||^2 = \sum_{\alpha} |\langle x, \phi_n \rangle|^2$$ holds, it suffices to show that the Parsevall identity holds for a dense subset of $L^2([0,1] \times [0,1])$. We consider the family $C([0,1] \times [0,1])$, which is dense in $L^2([0,1] \times [0,1])$ by the compactness of $[0,1] \times [0,1]$. So let $f \in C([0,1] \times [0,1])$. Then,
			\begin{align*}
			||f||^2 &= \langle f, f \rangle \\
			&= \int_0^1 \int_0^1 |f(x,y)|^2 dx dy \\
			&= \int_0^1 \Big[\sum_{n \in \N} \Big|\int_0^1 f(x,y) f_n(x) dx \Big|^2\Big] dy && \text{by completeness of $\{f_n\}$} \\
			&= \sum_{n \in \N} \Big[\int_0^1 \Big|\int_0^1 f(x,y) f_n(x) dx \Big|^2 dy \Big] && \text{since the sum converges to $f_y$} \\
			&= \sum_{m, n \in \N} \Big|\int_0^1 \int_0^1 f(x,y) f_n(x) dx g_m(y) dy \Big|^2 && \text{by completeness of $\{g_m\}$} \\
			&= \sum_{m, n \in \N} |\langle f, f_n g_m \rangle|^2,
			\end{align*} which shows that $\{f_n g_m\}_{n, m \in \N}$ is complete in $L^2([0,1] \times [0,1])$. \hfill $\Box$
			
			\par Note that continuity of $f$ and the compactness of $[0,1] \times [0,1]$ was tacitly used to guarantee the existence of all integrals that appear in the calculation. 
			
			\medspace
			
			Since the trigonometric family of polynomials of period $1$ and Bessel functions of the $0^{\text{th}}$ kind are each complete in $L^2([0,1])$, the family of functions given by their products is dense in $L^2([0,1] \times [0,1])$. Hence, we can represent $f(x,y,t)$ and any solution $u(x,y,t)$ by 
			\begin{align}\label{nonhom-f-rep}
			f(x,y,t) &= \sum_{m, n \in \N} F_{m,n}(t) J_0(\gamma_m x) \sin(n \pi y), \\ 
			u(x,y,t) &= \sum_{m, n \in \N} T_{m,n}(t) J_0(\gamma_m x) \sin(n \pi y). \label{non-hom-rep-u}
			\end{align} Here, note that the $F_{m,n}$ are determined and
			\begin{equation}\label{fmn}
			F_{m,n}(t) = \frac{1}{||J_0(\gamma m x) \sin(n \pi y) ||_2^2}\int\limits_0^1 \int\limits_0^1 x f(x,y,t) J_0(\gamma_m x) \sin(n \pi y) dx dy.
			\end{equation}
			Now, plugging into the differential equation for the nonhomogenous case yields%
			\begin{equation}\label{nonhom-induced} \left(\,\capd - \frac{\partial^2}{\partial x^2} - \frac{1}{x}  \frac{\partial}{\partial x} - \frac{\partial^2}{\partial y ^2} \right) \left[\sum_{m, n \in \N} T_{m,n}(t) J_0(\gamma_m x) \sin(n \pi y) \right] = \sum_{m, n \in \N} F_{m,n}(t) J_0(\gamma_m x) \sin(n \pi y).
			\end{equation} 
			By the computations done in the homogeneous case, we see that
			\begin{equation}\label{nonhom-helper}\left(\frac{\partial^2}{\partial x^2} + \frac{1}{x}  \frac{\partial}{\partial x} + \frac{\partial^2}{\partial y ^2}\right)\left[\sum_{m, n \in \N} T_{m,n}(t) J_0(\gamma_m x) \sin(n \pi y) \right] = -\sum_{m,n \in \N} T_{m,n}(t) ((n \pi)^2 + \gamma_m^2) J_0(\gamma_m x) \sin(n \pi y).
			\end{equation} 
			Combining (\ref{nonhom-induced}) and (\ref{nonhom-helper}) and noting that $\mu_{m,n} = (n \pi)^2 + \gamma_m ^2$, we obtain that
			$$\sum_{m, n \in \N} \Big(\,\capd T_{m,n}(t) + \mu_{m,n}T_{m,n}(t)\Big) J_0(\gamma_m x) \sin(n \pi y) =  \sum_{m, n \in \N} F_{m,n}(t) J_0(\gamma_m x) \sin(n \pi y).$$ Now, for any $l, k \in \N$, taking inner products with $J_0(\gamma_l x) \sin(k \pi y)$ and using the convergence of sums and linearity of the inner product shows that 
			\begin{align*}
			\sum_{m, n \in \N} \Big(\,\capd T_{m,n}(t) + \mu_{m,n}T_{m,n}(t)\Big) \langle J_0(\gamma_m x) \sin(n \pi y), J_0(\gamma_l x) \sin(k \pi y) \rangle \\
			=  \sum_{m, n \in \N} F_{m,n}(t) \langle J_0(\gamma_m x) \sin(n \pi y), J_0(\gamma_l x) \sin(k \pi y) \rangle.
			\end{align*} 
			Using the orthogonality of $\{J_0(\gamma_l x) \sin(k \pi x)\}$ shows that for all $m, n \in \N$,
			\begin{equation*}
			\capd T_{m,n}(t) + \mu_{m,n}T_{m,n}(t) = F_{m,n}(t),
			\end{equation*} 
			which is the Cauchy problem described in section 2.2 with $\lambda = \mu_{m,n}$ and $f(t) = F_{m,n}(t)$. Hence, we see that 
			$$T_{m,n}(t) =  b_{m,n}E_{\alpha, 1}(-\mu_{m,n} t^{\alpha})+c_{m,n}tE_{\alpha, 2}(-\mu_{m,n} t^{\alpha})+ \int_0^t (t-\xi)^{\alpha-1}E_{\alpha, \alpha}[-\mu_{m,n}(t-\xi)^{\alpha}]F_{m,n}(\xi)d\xi$$ 
			where $b_{m,n}$ are constant.
			\par We now determine these coefficients. Since $u(x,y,0) = 0$ and $$  \int_0^0 (t-\xi)^{\alpha-1}E_{\alpha, \alpha}[-\mu_{m,n}(t-\xi)^{\alpha}]F_{m,n}(\xi)d\xi = 0,$$ (\ref{non-hom-rep-u}) and (\ref{E-alpha-one-zero}) give that
			\begin{align*}
			&0 = u(x,y,0) = \sum_{m, n \in \N} \Big[b_{m,n}E_{\alpha, 1}(0) + 0c_{m,n}E_{\alpha,2}(0) \Big]J_0(\gamma_m x) \sin(n \pi y) \\
			&= \sum_{m,n \in \N} b_{m,n}J_0(\gamma_m x) \sin(n \pi y).
			\end{align*} Since the family of functions $\{J_0(\gamma_m x) \sin(n \pi y)\}$ are orthogonal, this shows that the $b_{m,n}$ must be identically $0$. To handle the values of $c_{m,n}$, first note that
			\begin{align*}
			&\frac{d}{dt}\Bigg( \int_0^t (t-\xi)^{\alpha-1} E_{\alpha, \alpha}[-\mu_{m,n}(t-\xi)^{\alpha}]F_{m,n}(\xi)d\xi\Bigg)_{t=0} \\ 
			&= (t - t)^{\alpha - 1}E_{\alpha, \alpha}(-\mu_{m,n}(t -t)^\alpha)F_{m,n}(t)\Big|_{t=0} + \Bigg(\int_0^t \frac{\partial}{\partial t} (t-\xi)^{\alpha -1}E_{\alpha, \alpha}(-\mu_{m,n}(t - \xi)^{\alpha})F(\xi)d\xi\Bigg)_{t=0} = 0.
			\end{align*} Then, the fact that $u_t(x, y, 0) = 0$ combined with (\ref{E-alpha-two-zero}), (\ref{nonhom-induced}), and identity (\ref{MittagDeriv}) gives that
			\begin{align*}
			0 &= u_t(x,y,0) \\
			&= \sum_{m,n}\Big[b_{m,n}\frac{d}{dt} E_{\alpha, 1}(-\mu_{m,n} t^{\alpha})+ c_{m,n}\frac{d}{dt}tE_{\alpha, 2}(-\mu_{m,n} t^{\alpha}) \Big]_{t=0} J_0(\gamma_m x) \sin(n \pi y) \\ 
			&= \sum_{m,n} \Big[b_{m,n} (-\mu_{m,n} t^{\alpha -1}) E_{\alpha, \alpha}(-\mu_{m,n} t^\alpha) + c_{m,n}\Big(t \frac{d}{dt} E_{\alpha, 2}(-\mu_{m,n} t^\alpha) + E_{\alpha, 2}(-\mu_{m,n} t^\alpha) \Big) \Big]_{t=0} J_0(\gamma_m x) \sin(n \pi y) \\
			&= \sum_{m,n} c_{m,n} J_0(\gamma_m x) \sin(n \pi y).
			\end{align*}
			As before, since $\{J_0(\gamma_m x)\sin(n \pi y)\}$ is an orthogonal family of functions, the constants $c_{m,n}$ must be identically $0$. This now gives us a final solution to the nonhomogeneous problem
			$$u(x,y,t) = \sum_{m,n \in \N} \left(\int_0^t (t-\xi)^{\alpha-1}E_{\alpha, \alpha}[-\mu_{m,n}(t-\xi)^{\alpha}]F_{m,n}(\xi)d\xi\right)J_0(\gamma_m x) \sin(n \pi x).$$ 
			Combining the solutions to the homogeneous and nonhomogeneous problems gives that a formal solution to the original problem given by (\ref{original-conditions-start}) - (\ref{original-conditions-end}),
			\begin{align}\label{full-solution}
			u(x,y,t) &= \sum_{m,n \in \N}T_{m,n}(t)J_0(\gamma_m x) \sin(n \pi y), \\
			T_{m,n}(t) &= C_{m,n} E_{\alpha, 1}(-\mu_{m,n} t^\alpha) + D_{m,n}t E_{\alpha, 2}(-\mu_{m,n} t^\alpha) + \int_0^t (t-\xi)^{\alpha-1}E_{\alpha, \alpha}[-\mu_{m,n}(t-\xi)^{\alpha}]F_{m,n}(\xi)d\xi,\label{full-solution-end}
			\end{align} where $C_{m,n},$ $D_{m,n},$ and $F_{m,n}$ are given by (\ref{cmn}), (\ref{dmn}), and (\ref{nonhom-f-rep}), respectively.
			
			\medskip
			
			
			\textbf{4. Uniqueness}
			
			We claim that any continuous solution to the problem (\ref{eq1}), (\ref{original-conditions-start}) - (\ref{original-conditions-end}) must be unique. To see this, suppose that $u(x,y,t)$ and $u'(x,y,t)$ solved this problem. Then, their difference $s(x,y,t) = u(x,y,t) - u'(x,y,t)$ would be a solution to
			\begin{equation*}
			\capd s(x,y,t) = s_{xx} + \frac{1}{x} s_x + s_{yy}
			\end{equation*}
			with boundary and initial conditions
			\begin{align*}
			\lim_{x \rightarrow 0} x s_x &= 0 &&s(1,y,t) = 0 \\
			s(x, 0, t) &= 0  &&s(x, 1, t) = 0 \\
			s(x, y, 0) &= 0   &&s_t(x,y,0) = 0.
			\end{align*} We are now in the same situation as we were in the nonhomogeneous problem, with $f(x,y,t) = 0$. The argument we used in the nonhomogeneous case now shows that
			\begin{align*}
			s(x,y,t) &= \sum_{m,n \in \N} \Big(\int_0^t (t-\xi)^{\alpha-1}E_{\alpha, \alpha}[-\mu_{m,n}(t-\xi)^{\alpha}]0d\xi\Big)J_0(\gamma_m x) \sin(n \pi x) \\
			&= 0.
			\end{align*} Thus $u = u'$.
			
			\medskip
			
			\textbf{5. Existence}
			\par We turn our attention now to showing that the formal solution given by (\ref{full-solution}) and (\ref{full-solution-end}) is well-defined and converges under proper assumptions on $f$, $\varphi$ and $\psi$. In particular, we will demonstrate the following result. 
			
			\textbf{Theorem.}\emph{ The formal solution given by (\ref{full-solution})-(\ref{full-solution-end}) to the problem given by (\ref{eq1})-(\ref{original-conditions-end}) is well-defined provided that:
				\begin{enumerate}
					\item $f_{tt}$ is continuous and $f(x,y,0) = 0$; 
					\item $f_{yyyxx}$ is bounded and for all $(x,y,t) \in \Omega$, $$f(x,0,t) = f(x,1,t) = f_{yy}(x,0,t) = f_{yy}(x,1,t) = f_{yyy}(0, y, t) = f_{yyy}(1,y,t) = f_{yyyx}(0,y,t) = 0;$$
					\item $f_{yxxxx}$ is bounded and for all $(x,y,t) \in \Omega$, $$f(x,0,t) = f(x,1,t) = f_{y}(0,y,t) = f_{y}(1,y,t) = f_{yx}(0, y, t) = f_{yxx}(0,y,t) = f_{yxx}(1,y,t) = f_{yxxx}(0,y,t) = 0;$$
					\item $\varphi_{yyyxx}$ is bounded and for all $x, y \in [0,1]$,  $$\varphi(x,0) = \varphi(x, 1) = \varphi_{yy}(x,0) = \varphi_{yy}(x, 1) = \varphi_{yyy}(0,y) = \varphi_{yyy}(1,y) = \varphi_{yyyx}(0,y) = 0;$$
					\item $\psi$ satisfies the same assumptions as $\varphi.$
				\end{enumerate}} The order of partial derivatives that $\varphi$ and $\psi$ must have exist could be changed, but doing so would require the boundary conditions they satisfy to change similarly. For the sake of simplicity, we leave them in this single case. Note also that five times continuous differentiability on $\varphi$ and $\psi$ implies the conditions demanded here.
				
				\par Section $3$ demonstrates that (\ref{full-solution})-(\ref{full-solution-end}) is a solution to (\ref{eq1})-(\ref{original-conditions-end}) if the differential operator $$\capd - \frac{1}{x}\frac{\partial}{\partial x} - \frac{\partial^2}{\partial x^2} - \frac{\partial^2}{\partial y^2}$$ can be applied term-by-term. To show this, it suffices to show that the series representing $u(x,y,t)$ converges uniformly on $[0,T]$ and that the formal series given by $u_{xx}$, $1/x u_x,$ $u_{yy},$ and $u_{tt}$ all converge uniformly on $[\epsilon, T]$ for all $\epsilon > 0$. The remainder of this section will be devoted to demonstrating this by applications of the Weierstrass M-test.
				
				\par We begin by giving a handful of useful simple facts that will be used throughout the bounding procedure.
				
				
				\textbf{Lemma 2.} If $f: D \rightarrow \mathbb{C}$ is a continuous function with Fourier-Bessel series $$\sum_{m,n \in \N} f_{m,n}(t)J_0(\gamma_m x) \sin(n \pi y),$$ then $$|f_{m,n}(t)| \leq ||f||_{\infty}$$ for all $m, n$, and $t$.
				
				\medskip
				
				\emph{Proof.} With $f$ and $f_{m,n}$ as given, $||f||_\infty$ is well-defined and finite by the compactness of $D$ and continuity of $f$. Now, for any $t$,
				\begin{align*}|f_{m,n}(t)|^2 &\leq \sum_{m,n \in \N} |f_{m,n}(t)|^2 && \text{and by Parseval,} \\
				&= \int_0^1 \int_0^1 x|f(x,y,t)|^2 dx dy \\
				&\leq \int_0^1 \int_0^1 ||f||_\infty^2 dx dy \\
				&= ||f||_\infty^2.
				\end{align*} Hence, $|f_{m,n}(t)| \leq ||f||_\infty$. \hfill $\Box$ \medskip
				
				
				
				\textbf{Lemma 3.}\emph{ If $f: [0, 1] \rightarrow \mathbb{C}$ and $f' \in L^1([0,1])$, then:
					\begin{enumerate}
						\item \[\int\limits_0^1 f(x) \cos(n \pi x)dx = \frac{-1}{n \pi} \int\limits_0^1 f'(x) \sin(n \pi x) dx;\]
						\item and if $f(0) = f(1) = 0$, then \[\int\limits_0^1 f(x) \sin(n \pi x) dx = \frac{1}{n \pi} \int\limits_0^1 f'(x)\cos(n \pi x) dx. \]
					\end{enumerate}}
					
					Note here that if $f$ being continuously differentiable or having a bounded derivative is enough to satisfy the requirement that $f' \in L^1([0,1])$.
					
					\emph{Proof.} 
					This is a simple application of integration by parts. \hfill $\Box$ \medskip
					
					
					
					
					The analogous bounding feature for $J_0(x)$ is contained in \cite{Tolstov}, Ch. 8, Section 19, and we state without proof.
					
					\textbf{Lemma 4.} \emph{ If $f: [0, 1] \rightarrow \mathbb{C}$ is $2k$ differentiable for $k \in \mathbb{N}$, $f(0) = f'(0) = f(1) = \dots  = f^{(2k -2)}(0) = f^{(2k-2)}(1) = f^{(2k-1)}(0) = 0$, and $f^{2k}$ is bounded, then there is a constant $C > 0$ such that for all $m$, \[\int\limits_0^1 \sqrt{x}f(x) J_0(\gamma_m x) dx \leq \frac{C M }{(\gamma_m)^{2k + 1/2}},\] where $M = \max(f^{(2k)}(x))$ and $\gamma_m$ is the $m$-th positive $0$ of $J_0(x)$. }
					
					
					We next give some useful facts about $F_{m,n}$ and its derivatives. Note that from the representation of $F_{m,n}$ in (\ref{fmn}), 
					\begin{align*}F'_{m,n}(t) &= \frac{d}{dt} \left( \frac{1}{||J_0(\gamma_m x) \sin(n \pi y)||_2^2} \int\limits_0^1 \int\limits_0^1 x f(x,y,t) J_0(\gamma_m x) \sin(n \pi y) dx dy\right) \\
					&= \frac{1}{||J_0(\gamma_m x) \sin(n \pi y)||_2^2} \int\limits_0^1 \int\limits_0^1 x f_t(x,y,t) J_0(\gamma_m x) \sin(n \pi y) dx dy,
					\end{align*} which is the $(m,n)$-th Fourier coefficient of $f_t(x,y,t)$ with respect to the family $\{J_0(\gamma_m x ) \sin(n \pi y)\}$. Similarly, $F''_{m,n}$ is the $(m,n)$-th Fourier coefficient of $f_{tt}(x,y,t)$. The representations of $F_{m,n}$, $F'_{m,n}$ and $F''_{m,n}$ as integrals of $f$, $f_t$ and $f_{tt}$ against $x J_0(\gamma_m x) \sin(n \pi y)$ demonstrate that the assumption $f_{tt}(x,y,t)$ is continuous implies that all $F_{m,n}$, $F'_{m,n}$, and $F''_{m,n}$ are continuous.
					

					\par Our final step before directly handling the series at hand is to place a bound on $|C_{m,n}|$ and $|D_{m,n}|.$ Recall from (\ref{cmn}) that $$C_{m,n} = \frac{1}{||J_0(\gamma_m x) \sin(n \pi y)||_2^2} \int_0^1\int_0^1 x J_0(\gamma_m x) \sin(n \pi y) \varphi(x,y) dxdy $$ where \begin{equation}\label{cmn-norm-part}
					||J_0(\gamma_m x) \sin(n \pi y)||_2^2 = \int_0^1\int_0^1 x|J_0(\gamma_m x) \sin(n \pi y)|^2 dxdy = \int_0^1|sin(n \pi y)|^2 \left(\int_0^1 x J_0^2(\gamma_m x) dx\right) dy.\end{equation} We bound (\ref{cmn-norm-part}) below. This argument will take two steps.
					
					To bound the component involving $J_0$, note that as a consequence of \cite{Tolstov}, eq. 14.1, there is a $C_1 > 0$ such that $$\frac{C_1}{\gamma_m} \leq \int_0^1 xJ_0^2( \gamma_m x) dx$$ for all $m \in \mathbb{N}$. Next, we use the fact from \cite{Tolstov}, Ch. 8, Section 10 that $$C_2 := \inf_{n \in \N} \frac{n\gamma_1}{\gamma_n}$$ is a nonzero, positive number, and hence that $$\frac{C_2}{\gamma_1 m}  \leq \frac{1}{\gamma_m}$$ for all $m \in \mathbb{N}$. Thus, taking $C_3$ to be $(C_1C_2)/\gamma_1$, for any $m \in \mathbb{N}$, \begin{equation}\label{final-bessel-lower-bound} \frac{C_3}{m} \leq \int_0^1 x J_0^2(\gamma_m x) dx.
					\end{equation} This completes our bounding of the $x$-component of (\ref{cmn-norm-part}).
					
					To handle the $y$-component, note that $$C_4 := \int_0^1 |\sin(\pi y)|^2dy $$ is evidently a positive real number. Hence, 
					\begin{align} 
					\int_0^1 |\sin(n \pi y)|^2 dy &\geq \int_0^{\frac{1}{n}} |\sin(n \pi y)|^2 dy && \text{and putting $t = ny$} \\
					&=\int_0^1 \frac{1}{n} |\sin(\pi t)|^2 dt \\
					&= \frac{C_4}{n} \label{final-sin-bound}
					\end{align} Combining (\ref{cmn-norm-part}), (\ref{final-bessel-lower-bound}), and (\ref{final-sin-bound}) demonstrates%
					\begin{equation}|C_{m,n}| \leq C nm \int_0^1\int_0^1x J_0(\gamma_m x) \sin(n \pi y) \varphi(x,y) dxdy,\label{init-cmn-bound}\end{equation} where each $C$ is a positive constant independent of $m$ and $n$. An identical argument can be applied to $|D_{m,n}|$ and $|F_{m,n}|$ to create similar bounds.
					
					
					
					\par We now prove the uniform convergence of (\ref{full-solution})-(\ref{full-solution-end}). Since $\sin(n \pi y)$ and $J_0(\gamma_m x)$ are both strictly bounded by constants, we can place the immediate bound 
					
					\begin{equation}\label{bounding-starting-point}
					|T_{m,n}(t) J_0(\gamma_m x) \sin(n \pi y) | \leq C_5 |T_{m,n}(t)|
					\end{equation} on the modulus of the terms of (\ref{full-solution})-(\ref{full-solution-end}). Furthermore, since $E_{\alpha, 1}(-\mu_{m,n} t^\alpha)$ and $E_{\alpha, 2}(-\mu_{m,n} t^\alpha)$ are real for real $t$, decrease as $t$ increases, and $E_{\alpha, 1}(0) = E_{\alpha, 2}(0) = 1$, we see that for all $m, n$ and $t \in \R$, 
					$$ |E_{\alpha, 1}(- \mu_{m,n} t^\alpha)| \leq 1,$$ and $$|tE_{\alpha, 2}(-\mu_{m,n} t^\alpha)| \leq t \leq T,$$
					and since $T$ is a fixed constant, we see from (\ref{bounding-starting-point}) that 
					\begin{align} |T_{m,n}(t)| &\leq C_6 \left|C_{m,n} + D_{m,n} + \int_0^t (t - \xi)^{\alpha-1} E_{\alpha, \alpha}[-\mu_{m,n}(t - \xi)^{\alpha}] F_{m,n}(\xi)d\xi\right| \\
					&\leq C_6\left(|C_{m,n}| + |D_{m,n}| + \left|\int\limits_0^t (t - \xi)^{\alpha -1} E_{\alpha, \alpha}[- \mu_{m,n}(t - \xi)^\alpha] F_{m,n}(\xi)d\xi \right| \right).\label{bounding-start-of-progress}
					\end{align}
					where $C_6$ is a positive constant.
					
					We begin by bounding the contributions of $|C_{m,n}|$ and $|D_{m,n}|$ to (\ref{bounding-start-of-progress}). Using the assumptions that $\varphi_{yyy}$ is continuous and $\varphi(x,0) = \varphi(x, 1) = \varphi_{yy}(x,0) = \varphi_{yy}(x,1) = 0$, we can apply Lemma 3 three times with respect to $y$ to (\ref{init-cmn-bound}) and obtain 
					\begin{equation}\label{cmn-n-bounded}|C_{m,n}| \leq \frac{C_{7} m}{n^2 \pi^3} \left|\int_0^1 \cos(n \pi y) \left(\int_0^1 x \varphi_{yyy}(x,y) J_0(\gamma_m x)  dx \right) dy \right|.\end{equation} 
					Now, since $\varphi_{yyy}$ is boundedly twice-differentiable with respect to $x$, and $\varphi_{yyy}(0,y) = \varphi_{yyy}(1,y) = \varphi_{yyyx}(0,y)$, $\sqrt{x} \varphi_{yyy}$ satisfies the same assumptions. Hence, 
					\begin{equation}\label{cmn-m-helper} \left| \int\limits_0^1 x \varphi_{yyy} (x,y) J_0(\gamma_m x) dx \right| = \left| \int\limits_0^1 \sqrt{ x} (\sqrt x \, \varphi_{yyy} (x,y)) J_0(\gamma_m x) dx \right| \leq \frac{C_8}{(\gamma_m)^{5/2}} \max_{x, y \in [0,1]}\left(\sqrt x\,\varphi_{yyyxx}(x,y)\right) \leq \frac{C_9}{m^{5/2}} \end{equation} by Lemma $4$. Combining (\ref{cmn-m-helper}) and (\ref{cmn-n-bounded}) now gives that
					$$|C_{m,n}| \leq \frac{C_{10}}{n^2 m^{3/2}} \int\limits_0^1 |\cos(n \pi y)| dy = \frac{C_{11}}{n^2 m^{3/2}}.$$
					and since 
					$$\sum_{m,n \in \N} \frac{1}{n^2 m^\frac{3}{2}}$$ 
					converges by the $p$-test, $\sum_{m,n \in \N}|C_{m,n}|$ is dominated by an absolutely convergent numerical series and hence is uniformly convergent by the Weierstrass M-test. 
					The proof of the uniform convergence of $\sum_{m,n \in \N} |D_{m,n}|$ is identical, and uses the same assumptions on $\psi$. 
					
					
					\par We now handle the contribution of the integral part to (\ref{bounding-start-of-progress}), which is given by 
					\begin{equation}\label{integral-series} \sum_{m,n \in \N} \left| \int_0^t (t - \xi)^{\alpha-1} E_{\alpha, \alpha}[-\mu_{m,n}(t - \xi)^{\alpha}] F_{m,n}(\xi)d\xi\right| . \end{equation} 
					After using (\ref{MittagDeriv}) to integrate by parts, we rewrite each term of (\ref{integral-series}) as 
					\begin{equation}\label{u-integral-starting-point} \bigg| \int_0^t \frac{1}{\mu_{m,n}} E_{\alpha,1}[-\mu_{m,n}(t - \xi)^\alpha] F'_{m,n}(\xi) d\xi - \frac{F_{m,n}(t)}{\mu_{m,n}} \bigg| \leq \left| \int_0^t \frac{1}{\mu_{m,n}} E_{\alpha,1}[-\mu_{m,n}(t - \xi)^\alpha] F'_{m,n}(\xi) d\xi \right| + \left| \frac{F_{m,n}(t)}{\mu_{m,n}} \right|.
					\end{equation}
					\par We handle each term in (\ref{u-integral-starting-point}). For the second, Bessel's inequality $ab \leq a^2 + b^2$ demonstrates that
					\begin{equation}\label{u-integral-second-bound} \left|\frac{F_{m,n}(t)}{\mu_{m,n}}\right| \leq \frac{1}{\mu_{m,n}^2} + |F_{m,n}(t)|^2.\end{equation} 
					Using the Parseval's identity and the fact that $\mu_{m,n} \geq mn$ now shows that \begin{equation}\label{u-first-integral-term-final}\sum_{n,m \in \N} \left(\frac{1}{\mu_{m,n}^2} + |F_{m,n}(t)|^2 \right) \leq ||f||_2^2 + \sum_{m,n \in \mathbb{N}} \frac{1}{(mn)^2} < \infty,\end{equation} completing the argument for the second term. For the first term, lemma $2$ and the assumption that $f_t(x,y,t)$ is continuous allows us to bound $F'_{m,n}(\xi)$ by a constant independent of $m$ and $n$. Thus, for any $\epsilon > 0$,
					\begin{align*}
					\left|\int_0^t \frac{1}{\mu_{m,n}} E_{\alpha,1}[-\mu_{m,n}(t - \xi)^\alpha] F'_{m,n}(\xi) d\xi \right| &\leq  \frac{C_{12}}{\mu_{m,n}}\left| \int_0^t E_{\alpha,1}[-\mu_{m,n}(t - \xi)^\alpha] d\xi \right| \\
					& \leq \frac{C_{12}}{\mu_{m,n}} \left(\int_0^{t-\epsilon} \frac{1}{1 + \mu_{m,n}(t-\xi)^\alpha} d\xi + \int_{t-\epsilon}^t 1 \, d\xi \right) && \text{by (\ref{godsend-identity})} \\
					& \leq \frac{C_{12}}{\mu_{m,n}} \left(\frac{T}{\mu_{m,n} \epsilon^{\alpha}} + \epsilon \right).
					\end{align*}
					If we take $\epsilon = \mu_{m,n}^{-1/3}$, then this shows 
					\begin{align}\label{almost-done-u-bound-first-part}
					\left|\int_0^t \frac{1}{\mu_{m,n}} E_{\alpha,1}[-\mu_{m,n}(t - \xi)^\alpha] F'_{m,n}(\xi) d\xi \right| &\leq \frac{C_{12}}{\mu_{m,n}} \left( \frac{T}{(\mu_{m,n})^{1-\alpha/3}} + (\mu_{m,n})^{-1/3}\right).
					\end{align}
					As $1 < \alpha < 2$ and $\mu_{m,n} > 1$, we see $({\mu_{m,n})^{-(1-\alpha/3)}}$ is bounded above by $(\mu_{m,n})^{-1/3}$. Thus, (\ref{almost-done-u-bound-first-part}) shows that
					\begin{align*} 
					\left|\int_0^t \frac{1}{\mu_{m,n}} E_{\alpha,1}[-\mu_{m,n}(t - \xi)^\alpha] F'_{m,n}(\xi) d\xi \right| \leq \frac{C_{13}}{(\mu_{m,n})^{4/3}}.
					\end{align*}
					Since $4/3 > 1$, the $p$-series test demonstrates that \begin{equation}\label{u-integral-part-first-final}\sum_{m,n \in \mathbb{N}} \frac{1}{(\mu_{m,n})^{4/3}} \leq \sum_{m,n \in \N} \frac{1}{(mn)^{4/3}} < \infty.\end{equation} Hence, combining (\ref{u-first-integral-term-final}) and (\ref{u-integral-part-first-final}), we see that (\ref{u-integral-starting-point}) is dominated by an absolutely convergent series, thus is uniformly convergent. We have now shown that series given by (\ref{full-solution}) is uniformly convergent.
					
					\medskip
					
					\par We next handle the uniform convergence of $u_{tt}$ on $[\delta, T]$ for a fixed $\delta > 0$. First, we compute $u_{tt}(x,y,t)$. Using (\ref{MittagDeriv}) and (\ref{MittagDerivGeneral}), $$\frac{d^2}{dt^2} E_{\alpha,1}(-\mu_{m,n}t^\alpha) = -\mu_{m,n}\frac{d}{dt} t^{\alpha-1}E_{\alpha,\alpha}(-\mu_{m,n}t^\alpha) = -\mu_{m,n}t^{\alpha-2}E_{\alpha,\alpha-1}(-\mu_{m,n}t^\alpha) $$ 
					and 
					$$ \frac{d^2}{dt^2}tE_{\alpha,2}(-\mu_{m,n}t^\alpha) = \frac{d}{dt} E_{\alpha,1}(-\mu_{m,n}t^{\alpha-1})= -\mu_{m,n}t^{\alpha-1} E_{\alpha,\alpha}(-\mu_{m,n}t^\alpha) .$$
					Lastly, following the computation done from (\ref{caputoint})-(\ref{caputo-second-der-end}) and using the assumption that $f(0) = 0$, we see that
					\begin{align*} 
					\frac{d^2}{dt^2} \int_0^t (t-\xi)^{\alpha-1} E_{\alpha,\alpha}&(-\mu_{m,n}(t-\xi)^\alpha)F_{m,n}(\xi) d\xi \\
					& = t^{\alpha-2} E_{\alpha,\alpha-1}(-\mu_{m,n}t^\alpha)F_{m,n}(0) + \int_0^t \frac{\partial}{\partial t} [z^{\alpha-2} E_{\alpha,\alpha-1}(-\mu_{m,n}z^\alpha)F_{m,n}(t-z)]dz \\
					& = \int_0^t z^{\alpha-2}E_{\alpha,\alpha-1}(-\mu_{m,n}z^\alpha)F_{m,n}'(t-z) dz.
					\end{align*} 
					We have now shown that 
					\begin{align*} 
					u_{tt}(x,y,t) &=  \sum_{m,n \in \N} T''_{m,n}(t) J_{0}(\gamma_m x) \sin(n \pi y) \\
					T''_{m,n}(t)	&= -\mu_{m,n}C_{m,n} t^{\alpha - 2}E_{\alpha, \alpha -1}(-\mu_{m,n}t^\alpha) -\mu_{m,n}D_{m,n}t^{\alpha - 1}E_{\alpha,\alpha}(-\mu_{m,n}t^\alpha) \\
					& \hspace{2cm} + \int_0^t z^{\alpha - 2}E_{\alpha,\alpha -1}(-\mu_{m,n}z^\alpha)f'(t-z) dz. \nonumber
					\end{align*}
					Our concern is now to bound each term $|T''_{m,n}(t) J_0(\gamma_m x) \sin(n \pi y)|$ above. Since $J_{0}(\gamma_m x) \sin(n \pi y)$ is uniformly bounded by a constant, we immediately have that 
					\begin{align*} |T''_{m,n}(t)J_0(\gamma_m x) \sin(n \pi y)| \leq  C_{14} \bigg[\big| & \mu_{m,n}C_{m,n} t^{\alpha - 2}E_{\alpha, \alpha -1}(-\mu_{m,n}t^\alpha)\big| + \big|\mu_{m,n}D_{m,n}t^{\alpha - 1}E_{\alpha,\alpha}(-\mu_{m,n}t^\alpha) \big|
					\\ & + \bigg|\int_0^t z^{\alpha - 2}E_{\alpha,\alpha -1}(-\mu_{m,n}z^\alpha)f'(t-z) dz \bigg|  \bigg].
					\end{align*} 
					We thus prove the uniform convergence of $u_{tt}$ by bounding each of the following terms with absolutely convergent numerical series:
					\begin{equation}\label{part1} 
					\left|  \mu_{m,n}C_{m,n} t^{\alpha - 2}E_{\alpha, \alpha -1}(-\mu_{m,n}t^\alpha)\right|,
					\end{equation} 
					\begin{equation} \label{part2}
					\left| \mu_{m,n}D_{m,n}t^{\alpha - 1}E_{\alpha,\alpha}(-\mu_{m,n}t^\alpha) \right|,
					\end{equation} 
					\begin{equation} \label{part3}
					\left|\int_0^t z^{\alpha - 2}E_{\alpha,\alpha -1}(-\mu_{m,n}z^\alpha)f'(t-z) dz \right|. 
					\end{equation}

					First consider (\ref{part1}). Applying (\ref{godsend-identity}) and the facts that $1< \alpha <2$ and $t > \delta$ on $[\delta, T]$,
					\begin{align*} \left|  \mu_{m,n}C_{m,n} t^{\alpha - 2}E_{\alpha, \alpha -1}(-\mu_{m,n}t^\alpha)\right| 
					&= \mu_{m,n}t^{\alpha-2}|C_{m,n}||E_{\alpha,\alpha-1}(-\mu_{m,n}t^\alpha)| \\
					&\leq C_{15}|C_{m,n}|\frac{\mu_{m,n}\delta^{\alpha - 2} }{1+\mu_{m,n}t^\alpha}  \\
					& \leq C_{16}|C_{m,n}| \frac{\mu_{m,n} \delta^{\alpha - 2} }{\mu_{m,n}\delta^\alpha}\\
					& = C_{17}|C_{m,n}|,
					\end{align*} where $C_{17}$ is a constant independent of $m,n$. Due to our previous results, this shows that (\ref{part1}) is bounded above by an absolutely convergent numerical series. Analogously, 
					\begin{align*} \left|  \mu_{m,n}D_{m,n} t^{\alpha - 1}E_{\alpha, \alpha}(-\mu_{m,n}t^\alpha)\right| 
					&= \mu_{m,n}t^{\alpha-1}|D_{m,n}||E_{\alpha,\alpha}(-\mu_{m,n}t^\alpha)| \\
					&\leq C'_{15}|D_{m,n}|\frac{\mu_{m,n}t^{\alpha - 1} }{1+\mu_{m,n}t^\alpha}  \\
					& \leq C'_{16}|D_{m,n}| \frac{\mu_{m,n} t^{\alpha - 1} }{\mu_{m,n}t^\alpha}\\
					& = C'_{16}|D_{m,n}|t^{-1} \\
					& \leq C'_{16}|D_{m,n}|\delta^{-1} \\
					& = C'_{17}|D_{m,n}|,
					\end{align*}
					which proves that (\ref{part2}) is also bounded above by an absolutely convergent numerical series.
					
					\par Finally, we consider (\ref{part3}). Integrating by parts with $dv = z^{\alpha - 2} E_{\alpha, \alpha - 1}(-\mu_{m,n} z^\alpha)$ and using (\ref{MittagDerivGeneral}) and (\ref{godsend-identity}) we see that 
					\begin{align} 
					&\left| \int_0^t z^{\alpha - 2}E_{\alpha,\alpha -1}(-\mu_{m,n}z^\alpha)F'_{m,n}(t-z) dz \right|\\
					 & \leq  \left|t^{\alpha-1}E_{\alpha,\alpha}(-\mu_{m,n}t^\alpha)F''_{m,n}(0)\right| + \left|\int_0^t z^{\alpha-1}E_{\alpha,\alpha}(-\mu_{m,n}z^\alpha)F''_{m,n}(t-z) dz\right| \nonumber \\
					&\leq C_{18} \frac{|F''_{m,n}(0)|}{\mu_{m,n}} + \left|\int_0^t z^{\alpha-1}E_{\alpha,\alpha}(-\mu_{m,n}z^\alpha)F''_{m,n}(t-z) dz\right| \label{utt-integral-start}.  
					\end{align} 
					The first term in (\ref{utt-integral-start}) is satisfactorily bounded by the same argument used in (\ref{u-integral-second-bound})-(\ref{u-first-integral-term-final}). For the second term, begin by noting that (\ref{strong-E-bound}) shows 
					\begin{equation}\label{E-alpha-alpha-factor-pull}z^{\alpha -1} E_{\alpha, \alpha }(- \mu_{m,n} z^\alpha) = \frac{-1}{z \mu_{m,n}}(-\mu_{m,n} z^\alpha E_{\alpha, \alpha} (- \mu_{m,n} z^\alpha)) = \frac{-1}{\mu_{m,n}} z^{-1}E_{\alpha, 0}(-\mu_{m,n} z^{\alpha}).\end{equation} After using Lemma 2 to estimate $F_{m,n}$ by a uniform constant, the integral term of (\ref{utt-integral-start}) then becomes bounded by
					\begin{equation}\label{Ea0-over-z-integral-bound-start}
					\frac{1}{\mu_{m,n}} \int\limits_0^t \left| z^{-1}E_{\alpha, 0}(-\mu_{m,n}t^\alpha) \right| dz. 
					\end{equation}
					Now, since $1/\Gamma(0)$ is understood to be $0$, 
					\begin{equation*}
					z^{-1}E_{\alpha, 0}(-\mu_{m,n}z^\alpha) = z^{-1}\sum_{k=1}^\infty\frac{(-\mu_{m,n}z^\alpha)^k}{\Gamma(k\alpha)} = \sum_{k = 1}^\infty \frac{(-\mu_{m,n})^kz^{\alpha k - 1}}{\Gamma(k \alpha)},
					\end{equation*}
					and since $\alpha > 1$, this demonstrates that $$\lim_{z \rightarrow 0} z^{-1} E_{\alpha, 0}(-\mu_{m,n}z^\alpha) = 0.$$ This demonstrates that $z^{-1}E_{\alpha, 0}(-\mu_{m,n}z^\alpha)$ can be continuously extended to the entirety of $[0,t]$. Furthermore, since $E_{\alpha, 0}$
					is diminishing as its argument decreases, for all $n, m \in \N$ and $z \in [0,t]$, 
					\begin{equation*}
					z^{-1} E_{\alpha, 0}(-\mu_{m,n}z^\alpha) \leq ||z^{-1}E_{\alpha, 0}(-\mu_{1,1}z^{\alpha})||_\infty,\end{equation*}
					which is finite by the compactness of $[0,t]$. Thus, for any $\epsilon \in (0, 1]$, applying (\ref{godsend-identity}) shows
					\begin{align*}
					\int_0^t\big|z^{-1} E_{\alpha, 0}(-\mu_{m,n}z^\alpha)\big| dz &= \int_0^{\epsilon}\big|z^{-1} E_{\alpha, 0}(-\mu_{m,n}z^\alpha)\big| dz + \int_{\epsilon}^t\big|z^{-1} E_{\alpha, 0}(-\mu_{m,n}z^\alpha)\big| dz \\
					&\leq \int_0^{\epsilon} C_{19} dz + \int_{\epsilon}^t \bigg| \frac{C_{20}}{\epsilon(1 + \mu_{m,n}z^\alpha)}\bigg| dz \\
					&\leq \epsilon C_{19} + \int_{\epsilon}^{t}\frac{C_{20}}{\mu_{m,n}\epsilon^{\alpha+1}} dz   \\
					&\leq \epsilon C_{19} + \frac{C_{21}}{\mu_{m,n}\epsilon^{\alpha + 1}}  && \text{and since $1 < \alpha < 2$,}\\
					&\leq \epsilon C_{19} + \frac{C_{21}}{\mu_{m,n}\epsilon^3}.
					\end{align*} 
					Taking $\epsilon = (\mu_{m,n})^{-1/4}$ now shows
					\begin{equation*}
					\left| \int\limits_0^t z^{\alpha - 1}E_{\alpha, \alpha}(-\mu_{m,n} t^{\alpha}) F''_{m,n}(t - z) dz \right| \leq \frac{C_{22}}{(\mu_{m,n})^{5/4}},
					\end{equation*} finishing the last portion of the bounding of (\ref{part3}).
					This completes the proof of the uniform convergence of $u_{tt}$.
					
					\medskip
					
					\par We next show that the series $u_{yy}(x,y,t)$ is bounded above termwise by an absolutely convergent numerical series on $[\delta,T]$ for an arbitrary $\delta > 0$.
					Taking the derivative twice with respect to $y$ shows that  
					\begin{align*} 
					u_{yy}(x,y,t) =   \sum_{m,n \in \N} (n \pi)^2 \Bigg(&C_{m,n}E_{\alpha,1}(-\mu_{m,n}t^\alpha) + D_{m,n}tE_{\alpha,2}(-\mu_{m,n}t^\alpha) \\
					& + \int_0^t (t-\xi)^{\alpha-1}E_{\alpha, \alpha}[-\mu_{m,n}(t-\xi)^{\alpha}]F_{m,n}(\xi)d\xi \Bigg) J_{0}(\gamma_m x) \sin(n \pi y).
					\end{align*} By placing a constant estimate on $J_0(\gamma_m x) \sin(n \pi y)$, we obtain that each term is bounded above by
					\begin{equation*}
					(n \pi)^2 \left|C_{m,n}E_{\alpha, 1}(- \mu_{m,n}t^\alpha))\right| + \left|D_{m,n} t E_{\alpha, 2}(-\mu_{m,n}t^\alpha)\right| + \left|\int\limits_0^t (t - \xi)^{\alpha - 1}E_{\alpha, \alpha}[-\mu_{m,n} (t - \xi)^\alpha]F_{m,n}(\xi) d\xi\right|.
					\end{equation*}
					As before, we bound this expression in three parts:
					\begin{equation} \label{part1-2}
					(n \pi)^2 \left|C_{m,n}E_{\alpha,1}(-\mu_{m,n}t^\alpha) \right|,
					\end{equation} 
					\begin{equation} \label{part2-2}
					(n \pi)^2 \left| D_{m,n}tE_{\alpha,2}(-\mu_{m,n}t^\alpha) \right|,
					\end{equation}
					\begin{equation}  \label{part3-2}
					(n \pi)^2 \left| \int_0^t (t-\xi)^{\alpha-1}E_{\alpha,\alpha}(-\mu_{m,n}(t-\xi)^\alpha)F_{m,n}(\xi) d\xi \right|.
					\end{equation}
					
					For (\ref{part1-2}), note that since we only consider $t \in [\delta, T]$ an application of (\ref{godsend-identity}) demonstrates 
					\begin{align} \label{cmn-second-der-bd} 
					(n \pi)^2 \left|C_{m,n}E_{\alpha,1}(-\mu_{m,n}t^\alpha) \right| \leq (n \pi)^2 \frac{C_{23}}{\mu_{m,n}\delta^\alpha } \leq C_{24}|C_{m,n}|,
					\end{align} as $\mu_{m,n} \geq n^2$ for all $n$. For (\ref{part2-2}), an analogous argument shows that \begin{equation}\label{dmn-second-der-bd}(n \pi)^2\left| D_{m,n}tE_{\alpha,2}(-\mu_{m,n}t^\alpha) \right| \leq C_{25}|D_{m,n}|.\end{equation}
					As we have previously sufficiently bounded $|C_{m,n}|$ and $|D_{m,n}|$ in the proof of the convergence of $u$, this finishes everything for $u_{yy}$ except (\ref{part3-2}). 
					
					For (\ref{part3-2}), notice that as in (\ref{E-alpha-alpha-factor-pull}), $$(t - \xi)^{\alpha -1} E_{\alpha, \alpha}(- \mu_{m,n} (t - \xi)^\alpha) = \frac{1}{\mu_{m,n}} (t - \xi)^{-1} E_{\alpha, 0}(-\mu_{m,n} (t - \xi )^\alpha).$$ This demonstrates that
					\begin{align*}
					(n \pi)^2 \left| \int_0^t (t-\xi)^{\alpha-1}E_{\alpha,\alpha}(-\mu_{m,n}(t-\xi)^\alpha)F_{m,n}(\xi) d\xi \right| &= \frac{(n \pi)^2}{\mu_{m,n}} \int\limits_0^t \left|(t - \xi)^{-1} E_{\alpha, 0} (-\mu_{m,n}(t - \xi)^\alpha)\right| \, \left|F_{m,n}(\xi)\right| d\xi \\
					&\leq ||F_{m,n}||_\infty \frac{(n \pi)^2}{\mu_{m,n}} \int\limits_0^t\left|(t - \xi)^{-1} E_{\alpha, 0} (-\mu_{m,n}(t - \xi)^\alpha)\right| d\xi. 
					\end{align*} We have now arrived at exactly the integral that appears in (\ref{Ea0-over-z-integral-bound-start}), and hence can follow the same argument to have its value bounded by $C_{25} (\mu_{m,n})^{-1/4}$, leaving us with
					\begin{equation}\label{sad-times}
					(n \pi)^2 \left| \int_0^t (t-\xi)^{\alpha-1}E_{\alpha,\alpha}(-\mu_{m,n}(t-\xi)^\alpha)F_{m,n}(\xi) d\xi \right| \leq ||F_{m,n}||_\infty \frac{(n \pi)^2}{(\mu_{m,n})^{5/4}}.
					\end{equation} We now apply the conditions we placed on $f$ to eliminate the $(n \pi)^2$ term. Throughout the following argument, let $\xi \in [\delta, T]$ be given. Note that at the start of this section, $F_{m,n}(\xi)$ was bounded by \begin{equation}\label{fmn-ugly-bound}|F_{m,n}(\xi)| \leq C_{26} m n \int\limits_0^1 \int\limits_0^1 x f(x,y,\xi) J_0(\gamma_m x) \sin(n \pi y) dx dy. \end{equation} Since $f_{yyy}$ is assumed to exist, and $f(x,0,t) = f(x,1,t) = f_{yy}(x,0,t) = f_{yy}(x,1,t) = 0$, applying lemma $3$ three times yields that 
					\begin{equation}\label{fmn-bound-start} |F_{m,n} (\xi)| \leq \frac{C_{26}m}{n^2 \pi^3} \int\limits_0^1 \cos(n \pi y) \left( \int\limits_0^1 x f_{yyy}(x,y,\xi) J_0(\gamma_m x) dx\right) dy. \end{equation} Now, as $f_{yyyxx}$ exists, is bounded and satisfies $f_{yyy}(0, y, t) = f_{yyy}(1, y, t) = f_{yyyx}(0, y, t) = 0,$ the same holds of $\sqrt{x} f_{yyy}$. Hence, by lemma $4$, and discussion of the zeros of $J_0(x)$, \begin{equation}\label{fmn-m-helper} \left| \int\limits_0^1 x f_{yyy} (x,y, \xi) J_0(\gamma_m x) dx \right| = \left| \int\limits_0^1 \sqrt{ x} (\sqrt x \, f_{yyy} (x,y\xi)) J_0(\gamma_m x) dx \right| \leq \frac{C_{27}}{(\gamma_m)^{5/2}} \max(\sqrt{x}f_{yyyxx}) \leq \frac{C_{28}}{m^{5/2}}. \end{equation} Bringing together (\ref{fmn-bound-start}) and (\ref{fmn-m-helper}) thus demonstrates that $$|F_{m,n}(\xi)| \leq \frac{C_{29}}{n^2 m^{3/2}},$$ and since $\xi$ was arbitrary, $||F_{m,n}||_\infty$ has the same bound. Therefore,
					$$ 	(n \pi)^2 \left| \int_0^t (t-\xi)^{\alpha-1}E_{\alpha,\alpha}(-\mu_{m,n}(t-\xi)^\alpha)F_{m,n}(\xi) d\xi \right| \leq \frac{C_{29}}{(\mu_{m,n})^{5/4}m^{3/2}},$$ which is a sufficient bound to complete the bounding of (\ref{part3-2}) and thus the uniform convergence of $u_{yy}$.
					\medskip

					
					\par Next, we show the uniform convergence of the series that represents a term-by-term application of $\frac{1}{x}\frac{\partial}{\partial x}$ to $u$ for $t$ on the interval $[\delta,T]$:
					\begin{equation*} 
					\frac{1}{x}u_x(x,y,t) = \sum_{m,n \in \N} T_{m,n}(t)\left(-\gamma_m \frac{1}{x}J_1(\gamma_m x)\right) \sin(n \pi y).
					\end{equation*} 
					Note that (\ref{Bessel-derivative-pos-shift}) was used to compute the derivative of $J_0{\gamma_m x}$ here. Using property (\ref{bessel-division-killer}) of the Bessel function, we rewrite the sum as
					\begin{equation*} 
					\sum_{m,n \in \N} T_{m,n}(t)\left( \frac{-\gamma_m }{2}(J_0(\gamma_m x) + J_2(\gamma_m x)) \right)\sin(n \pi y).
					\end{equation*} Once we estimate away $J_0(\gamma_m x), J_2(\gamma_m x)$, and $\sin(n \pi x)$ with a uniform constant, the terms that we must bound above in order to apply the Weierstrass M-test are
					\begin{equation}
					\frac{\gamma_m}{2} \big|  C_{m,n}E_{\alpha,1}(-\mu_{m,n}t^\alpha) \big|, \label{part1-3}
					\end{equation}
					\begin{equation}
					\frac{\gamma_m}{2} \big| D_{m,n}tE_{\alpha,2}(-\mu_{m,n}t^\alpha) \gamma_m \big|, \label{part2-3} \end{equation}
					\begin{equation}\frac{\gamma_m}{2} \left| \int_0^t (t-\xi)^{\alpha -1} E_{\alpha,\alpha}(-\mu_{m,n}(t-\xi)^\alpha)F_{m,n}(\xi) d\xi \right|.\label{part3-3}
					\end{equation}
					
					For (\ref{part1-3}) and (\ref{part2-3}), an analogous argument to as in (\ref{cmn-second-der-bd})-(\ref{dmn-second-der-bd}) demonstrates that 
					\begin{align*}
					\frac{\gamma_m}{2} \big|C_{m,n}E_{\alpha,1}(-\mu_{m,n}t^\alpha) \big| &\leq C_{30}|C_{m,n}|, \\
					\frac{\gamma_m}{2} \big|D_{m,n}tE_{\alpha,2}(-\mu_{m,n}t^\alpha) \gamma_m \big| &\leq C_{31} |D_{m,n}|.
					\end{align*} Now, since $|C_{m,n}|$ and $|D_{m,n}|$ have been bounded above by terms in absolutely convergent series in the proof of $u$'s convergence, this fulfills our obligation for (\ref{part1-3}) and (\ref{part2-3}).
					\par Last, we must bound (\ref{part3-3}). As the integral is the same as in (\ref{part3-2}), we are able to obtain the bound  
					\begin{equation}\label{sad-times-2}
					\frac{\gamma_m}{2} \left| \int_0^t (t-\xi)^{\alpha-1}E_{\alpha,\alpha}(-\mu_{m,n}(t-\xi)^\alpha)F_{m,n}(\xi) d\xi \right| \leq ||F_{m,n}||_\infty \frac{\gamma_m}{2(\mu_{m,n})^{5/4}} \leq C_{32}||F_{m,n}||_\infty \frac{m}{(\mu_{m,n})^{5/4}}.
					\end{equation}
					
					We now leverage our assumptions on $f$ to cancel the $m$ term in (\ref{sad-times-2}). The process is similar to as in (\ref{sad-times}): for any $\xi$, since $$f(x,0,t) = f(x, 1,t) = f_{y}(0,y,t) = f_{y}(1, y, t) = f_{yx}(0, y, t) = 0$$ and $f_{yxx}$ is bounded, applying lemmas $3$ then $4$ to (\ref{fmn-ugly-bound}) shows 
					\begin{equation*}
					|F_{m,n}(\xi)| \leq \frac{C_{33}}{m^{3/2}},
					\end{equation*} 
					and hence that $||F_{m,n}||_\infty$ has the same bound. Thus,
					$$\frac{\gamma_m}{2} \left| \int_0^t (t-\xi)^{\alpha-1}E_{\alpha,\alpha}(-\mu_{m,n}(t-\xi)^\alpha)F_{m,n}(\xi) d\xi \right| \leq C_{34} \frac{1}{(\mu_{m,n})^{5/4}m^{1/2}},$$ which is a strong enough bound to guarantee the uniform convergence of $(1/x)u_x$. 
					
					\medskip
					
					\par To finish the existence portion, we show the uniform convergence of $u_{xx}(x,y,t)$ on $[\delta,T]$. First, we compute $(J_0(\gamma_m x))_{xx}$. By applying (\ref{Bessel-derivative-pos-shift}) and then (\ref{bessel-derivative-fact}), we see that
					\begin{align*}\frac{\partial^2}{\partial x^2} J_0(\gamma_m x) &= \frac{\partial}{\partial x}\left( \gamma_m J_1 (\gamma_m x) \right), \\
					&= \frac{-\gamma_m^2}{2} \left(J_0(\gamma_m x) - J_2(\gamma_m x)\right).
					\end{align*} This now shows that
					\begin{equation*} 
					u_{xx}(x,y,t) = \sum_{m,n \in \N} \frac{-\gamma_m^2}{2} T_{m,n}(t)\left(J_2(\gamma_m x) - J_0(\gamma_m x)\right) \sin(n \pi y).
					\end{equation*} 
					We again estimate  $J_0(\gamma_m x), J_2(\gamma_m x)$, and $\sin(n \pi x)$ with a uniform constant, to obtain that we must bound
					\begin{equation}
					\frac{\gamma_m^2}{2} \big|  C_{m,n}E_{\alpha,1}(-\mu_{m,n}t^\alpha) \big|, \label{part1-final}
					\end{equation}
					\begin{equation}
					\frac{\gamma_m^2}{2} \big| D_{m,n}tE_{\alpha,2}(-\mu_{m,n}t^\alpha) \gamma_m \big|, \label{part2-final} \end{equation}
					\begin{equation}\frac{\gamma_m^2}{2} \left| \int_0^t (t-\xi)^{\alpha -1} E_{\alpha,\alpha}(-\mu_{m,n}(t-\xi)^\alpha)F_{m,n}(\xi) d\xi \right|\label{part3-final}
					\end{equation} with absolutely convergent series.
					
					Once more, (\ref{part1-final}) and (\ref{part2-final}) are handled by an analogous argument to as in (\ref{cmn-second-der-bd})-(\ref{dmn-second-der-bd}) which demonstrates that 
					\begin{align*}
					\frac{\gamma_m^2}{2} \big|C_{m,n}E_{\alpha,1}(-\mu_{m,n}t^\alpha) \big| &\leq C_{35}|C_{m,n}|, \\
					\frac{\gamma_m^2}{2} \big|D_{m,n}tE_{\alpha,2}(-\mu_{m,n}t^\alpha) \gamma_m \big| &\leq C_{36} |D_{m,n}|.
					\end{align*} Now, since $|C_{m,n}|$ and $|D_{m,n}|$ have been bounded above by terms in absolutely convergent series in the proof of $u$'s convergence, this fulfills our obligation for (\ref{part1-3}) and (\ref{part2-3}).
					\par Finally, we must place a bound on (\ref{part3-final}). The integral in the term is the same as in (\ref{part3-2}), hence we are able to obtain the bound  
					\begin{equation*}
					\frac{\gamma_m^2}{2} \left| \int_0^t (t-\xi)^{\alpha-1}E_{\alpha,\alpha}(-\mu_{m,n}(t-\xi)^\alpha)F_{m,n}(\xi) d\xi \right| \leq ||F_{m,n}||_\infty \frac{\gamma_m^2}{2(\mu_{m,n})^{5/4}} \leq C_{37}||F_{m,n}||_\infty \frac{m^2}{(\mu_{m,n})^{5/4}}.
					\end{equation*} 
					We will now apply the remaining unused hypotheses on $f$ to cancel the $m^2$ term. Our method is similar to as in (\ref{sad-times}). If $\xi \in [\delta, T]$ is given, since $$f(x,0,t) = f(x, 1,t) = f_{y}(0,y,t) = f_{y}(1, y, t) = f_{yx}(0, y, t) = f_{yxx}(0, y, t) = f_{yxx}(1,y,t) = f_{yxxx}(0,y,t) = 0$$ and $f_{yxxxx}$ is bounded, applying lemmas $3$ then $4$ twice to (\ref{fmn-ugly-bound}) shows 
					\begin{equation*}
					|F_{m,n}(\xi)| \leq \frac{C_{33}}{m^{7/2}},
					\end{equation*} and hence that $||F_{m,n}||_\infty$ has the same bound. Thus,
					$$\frac{\gamma_m^2}{2} \left| \int_0^t (t-\xi)^{\alpha-1}E_{\alpha,\alpha}(-\mu_{m,n}(t-\xi)^\alpha)F_{m,n}(\xi) d\xi \right| \leq C_{34} \frac{1}{(\mu_{m,n})^{5/4}m^{3/2}},$$ which bounds (\ref{part3-final}) above with an absolutely convergent numerical series. This completes the convergence of $u_{xx}$ and the section on existence.
					\medskip
					
					
					\textbf{6. Stability}
					
					\par We now demonstrate that our solution depends continuously on the choice of initial conditions. To do this, it suffices to show an estimation of $||u(x,y,t)||$ in terms of $||\varphi||, ||\psi||,$ and $||f||$ with no additive constants. This is due to the linearity of the homogeneous and nonhomogeneous problems: if $u$ solves (\ref{eq1})-(\ref{original-conditions-end}) with $\varphi$, $\psi$, and $f$ given, and $u'$ solves (\ref{eq1})-(\ref{original-conditions-end}) with $\varphi'$, $\psi'$, and $f'$ given, then $u - u'$ is given by the solution to (\ref{eq1})-(\ref{original-conditions-end}) with $\varphi - \varphi'$, $\psi - \psi'$, and $f - f'$ given. Our approach is largely similar to the methods used in existence.
					
					\par Since any solution to the original problem (\ref{eq1})-(\ref{original-conditions-end}) is given by the sum of solutions to the homogeneous problem (\ref{hom-eq})-(\ref{hom-conditions-end}) and nonhomogeneous problem (\ref{nonhom-eq})-(\ref{nonhom-conditions-end}), we show bounds on each solution individually.
					
					
					\par The solution to the homogeneous problem (\ref{hom-eq})-(\ref{hom-conditions-end}) is given by $$ u(x,y,t) = \sum_{m, n \in \N} [C_{m,n} E_{\alpha, 1}(-\mu t^\alpha) + D_{m,n}t E_{\alpha, 2}(-\mu t^\alpha)]J_0(\gamma_m x) \sin(n \pi y),$$ where $C_{mn}$ and $D_{mn}$ are as in (\ref{cmn}) and (\ref{dmn}). Using Parseval's identity on the $x$ and $y$ components shows that
					\begin{align*}
					&||u(\cdot ,\cdot , t)||_2^2 = \sum_{m,n \in \N} |C_{m,n} E_{\alpha, 1}(-\mu_{m,n}t^\alpha) + D_{m,n}tE_{\alpha, 2}(-\mu_{m,n}t^\alpha)|^2 \\
					&\leq \sum_{m,n \in \N} \Big[|C_{m,n}E_{\alpha, 1}(-\mu_{m,n} t^\alpha)|^2 + 2|C_{m,n}E_{\alpha, 1}(-\mu_{m,n} t^\alpha)||D_{m,n}tE_{\alpha, 2}(-\mu_{m,n}t^\alpha)| + |D_{m,n}tE_{\alpha, 2}(-\mu_{m,n}t^\alpha)|^2\Big] \\
					&\leq \sum_{m,n \in \N}2\Big[|C_{m,n}E_{\alpha, 1}(-\mu_{m,n} t^\alpha)|^2 + |D_{m,n}tE_{\alpha, 2}(-\mu_{m,n}t^\alpha)|^2\Big],
					\end{align*} 
					where the last step is an application of Bessel's inequality. Now, since 
					$$E_{\alpha, 2}(- \mu_{m,n}t^\alpha) \leq \frac{1}{1 + \mu_{m,n}t^\alpha} \leq \frac{1}{1 + t^\alpha},$$
					and $t \geq 0$, we see that 
					$$[E_{\alpha, 2}(-\mu_{m,n}t^\alpha)]^2 \leq \frac{1}{(1+t^\alpha)^2} \leq \frac{1}{1 + t^{2\alpha}}.$$
					Hence, 
					$$|D_{m,n}t E_{\alpha, 2}(-\mu_{m,n}t^\alpha)|^2 \leq |D_{m,n}|^2 \frac{t^2}{1 + t^{2\alpha}}.$$
					In the case that $0\leq t \leq 1$, certainly $t^2\leq1$, and therefore
					$$|D_{m,n}t E_{\alpha, 2}(-\mu_{m,n}t^\alpha)|^2  \leq |D_{m,n}|^2.$$
					Similarly, when $t\geq 1 $ then $t^2\leq t^{2\alpha}$, and we again obtain the above estimate.
					Since the $E_{\alpha, 1}$ term can be bounded above by $1$, combining our estimates yields
					\begin{align*}
					||u(\cdot, \cdot, t)||_2^2 &\leq 2\left[\sum_{m,n \in \N}|C_{m,n}|^2 + \sum_{m,n \in \N}|D_{m,n}|^2\right] \\
					&= 2 ||\varphi||_2^2 + ||\psi||_2^2.
					\end{align*} Hence, $$||u(\cdot \, \cdot \, t)||_2 \leq \sqrt{2(||\varphi||_2^2 + ||\psi||_2^2)},$$ completing the homogeneous case.
					

					\par For the nonhomogeneous case, we have 
					$$u(x,y,t) = \sum_{m,n \in \N} \left(\int_0^t (t-\xi)^{\alpha-1}E_{\alpha, \alpha}[\mu_{m,n}(t-\xi)^{\alpha}]F_{m,n}(\xi)d\xi\right)J_0(\gamma_m x) \sin(n \pi x).$$ Estimating the $x$ and $y$ portions of $u$ by constants then gives $$|u(x,y,t)| \leq \sum_{m,n \in \N}\bigg|\int_{0}^t (t - \xi)^{\alpha - 1}E_{\alpha, \alpha}[-\mu_{m,n}(t-\xi)^\alpha]F_{m,n}(\xi) d\xi\bigg|.$$ 
					Now, for any $m$ and $n$, we can apply the identity (\ref{godsend-identity}) and then lemma $2$ to obtain
					\begin{align}
					\bigg|\int_0^t (t - \xi)^{\alpha - 1}E_{\alpha, \alpha}(-\mu_{m,n}(t- \xi)^\alpha)F_{m,n}(\xi)d\xi\bigg| &= \frac{1}{\mu_{m,n}}\bigg|\int_0^t(t - \xi)^{-1}E_{\alpha, 0}(-\mu_{m,n}(t-\xi)^\alpha) F_{m,n}(\xi)d\xi\bigg|\\
					&\leq \frac{||f||_\infty}{\mu_{m,n}}\int_0^t\bigg|(t - \xi)^{-1}E_{\alpha, 0}(-\mu_{m,n}(t-\xi)^\alpha)\bigg|d\xi \label{nonhom-stability-start}
					\end{align} 
					We are now in precisely the situation of (\ref{Ea0-over-z-integral-bound-start}), and following the same argument shows that for any $\epsilon \in (0, 1]$,
					\begin{equation*}
					\int_0^t\big|(t- \xi)^{-1} E_{\alpha, 0}(-\mu_{m,n}(t-\xi)^\alpha)\big| d\xi \leq \frac{C_1 T }{\mu_{m,n}\epsilon^{\alpha + 1}} + \epsilon C_2.
					\end{equation*}
					Here, the factor of $T$ demonstrates that the estimate that we obtain on $u$ will be dependent on the geometry of the domain. To obtain an estimate that depends minimally on the domain, we minimize this factor by taking 
					$$\epsilon = \left(\frac{T^*}{\mu_{m,n}}\right)^{{\frac{1}{4}}}$$ where $T^* = \max(T, 1)$. Then, 
					\begin{align*}
					\int_0^t\bigg|(t - \xi)^{-1}E_{\alpha, 0}(-\mu_{m,n}(t-\xi)^\alpha)\bigg|d\xi \leq \frac{(C_1 T^*)}{(\mu_{m,n})(\mu_{m,n}^{-\frac{3}{4}})((T^*)^{\frac{3}{4}})} + \frac{C_2 (T^*)^{\frac{1}{4}}}{(\mu_{m,n})^\frac{1}{4}} \leq \frac{C_3 (T^*)^{\frac{1}{4}}}{(\mu_{m,n})^{\frac{1}{4}}},
					\end{align*} which, combined with (\ref{nonhom-stability-start}) shows that $$|u(x,y,t)| \leq C_3||f||_\infty (T^*)^{\frac{1}{4}}\sum_{m,n \in \N} \frac{1}{(\mu_{m,n})^{\frac{5}{4}}},$$ and hence that $$||u||_\infty \leq C_4 ||f||_\infty (T^*)^{\frac{1}{4}},$$ completing the non-homogeneous case.
					
					\medskip
					
					{\bf Acknowledgement.} Authors are grateful to Dr Erkinjon Karimov for his useful remarks.\\
					This material is based upon work supported by the National Science Foundation under Grant No. NSF 1658672.
					
					\medskip
					
					\bibliography{References}

\textbf{References}
					
					\begin{enumerate}
						
						\bibitem{Sob74} Sobolev S.L. Introduction to the theory of cubature formulas. -Moscow: Nauka, 1974, 808 p.
						\bibitem{SobVas96} Sobolev S.L., Vaskevich V.L. Cubature formulas. -Novosibirsk:
						Publishing house of IM Siberian branch of RAS, 1996, -484 p.
						\bibitem{KilSrivTru} Kilbas A.A., Srivastava H.M. and Trujillo J.J. Theory and Applications of Fractional Differential Equations, North-Holland Mathematical Studies, Vol. 204, Elsevier, Amsterdam, 2006.
						\bibitem{Magin} R.L. Magin, Fractional Calculus in Bioengineering, Begell House Publisher, Inc. Connecticut,
						2006.
						\bibitem{WestBolGrig} West B.J.,  Bologna M. and Grigolini P. Physics of Fractal operators, New York, Springer, 2003.
						\bibitem{JesusMach} Jesus I.S., Tenreiro Machado J.A. Fractional control of heat diffusion systems, Nonlin. Dyn., 54(3), (2008), 263--282.
						\bibitem{MillerRoss} Miller K.S. and Ross B. An Introduction to the Fractional Calculus and Fractional Differential Equations, Wiley, New York 1993.
						\bibitem{Kilbas} Kilbas, A.A.  Boundary Value Problems, Special Functions and Fractional Calculus, Belorussian State University, Minsk 1996. (ISBN 985-6144-40- X) (Proc. Int. Conf., 90-th Birth Anniversary of Academician F.D. Gakhov, Minsk, Belorussia, 16-20 February 1996)
						\bibitem{Ross} Ross, B.  Fractional Calculus and its Applications, Lecture Notes in Mathematics 457, Springer Verlag, Berlin 1975. (Proc. Int. Conf. held at Univ. of New Haven, USA, 1974)
						\bibitem{Mainardi} Mainardi F. The fundamental solutions for the fractional diffusion-wave equation, Appl. Math. Lett. 9(6), (1996), 23--28.
						\bibitem{Podlubny} Podlubny I. Fractional Differential Equations, Academic Press, New York, 1999.
						\bibitem{OldSpan} Oldham K.B., Spanier J. The Fractional Calculus, Academic Press, New York,
						1974.
						\bibitem{BoyLuch} Boyadijev L., Luchko Y. Multi-dimensional $\alpha$-fractional diffusion-wave equation and some properties of its fundamental solution, Computers and Mathematics with Applications (2017)
						\bibitem{OzKar} Ozdemir N., Karadeniz D. Fractional diffusion-wave problem in cylindrical coordinates, Physics Letters A 372 (2008), http://dx.doi:10.1016/j.physleta.2008.07.054
						\bibitem{CiaRon} Ciaurri O., Roncal L. The wave equation for the Bessel Laplacian, J. Math. Anal. Appl. 409 (2014), http://dx.doi.org/10.1016/j.jmaa.2013.06.039
						\bibitem{BazBaz} Bazhlekova E., Bazhlekov I. Subordination approach to multi-term time-fractional diffusion-wave equations, Journal of Computational and Applied Mathematics (2017)
						\bibitem{Karimov} Karimov E., Pirnafasov S. Higher order multi-term time-fractional partial differential equations involving Caputo-Fabrizio derivative, Electronic Journal of Differential Equations, Vol. 2017 (2017), No. 243, pp. 1--11
						\bibitem{Pov} Povstenko Y. Neumann boundary-value problems for a time-fractional diffusion-wave equation in a half-plane, Computers and Mathematics with Applications (2012) 
						\bibitem{MartReySosa} Yepez-Martinez H., Reyes J.M., Sosa I.O. Analytical solutions to the fractional wave equation with variable dielectric function, Lat. Am. J. Phys. Educ. Vol. 8, No. 1, March (2014)
						\bibitem{Tolstov} Tolstov G.P. Fourier Series, Dover Publications Inc., New York, 1962
							\bibitem{Hilfer} R.Hilfer (Ed) Applications of Fractional Calculus in Physics, World Sci. Publishing, New York, 2000.
							\bibitem{Al-Musalhi} F.Al-Musalhi, N.Al-Salti, S.Kerbal. Inverse problems of a Fractional Differential Equation with Bessel Operator. Math. Model. Nat. Phenom.Vol.12, No 3, 2017, pp.105--113.
							\bibitem{Ruzhansky} P.Agarwal, E.Karimov, M.Mamchuev, M.Ruzhansky. On Boundary-Value Problems for a Partial Differential Equation with Caputo and Bessel Operators. Recent Applications of Harmonic Analysis to Functional Spaces, Differential Equations, and Data Science. Novel Methods in Harmonic Analysis, Vol. 2, 2017, pp.707-717
						
						\end{enumerate}

					%
					%
					
					\begin{tabular}{p{10cm}l}
				\textbf{Joseph David}; Florida State University; Tallahassee, FL, USA&\\
				E-mail: jkd14d@my.fsu.edu & \\
				\textbf{Alexander Nolte}; Tufts University, Medford, MA, USA & \\
				E-mail: alexander.nolte@tufts.edu &\\
				\textbf{Julie Sherman}; University of Minnesota - Twin Cities; Minneapolis, MN, USA;         & \\
						E-mail: sherm322@umn.edu     &  \\
					\end{tabular}
					%
					%
					\label{lastpage}

				\end{document}